\documentclass[11pt,twoside,a4paper]{article}

\usepackage{amsmath,bbm,amssymb,mathtools,geometry, xcolor,mathabx,mathrsfs,amsthm, hyperref}
\usepackage{verbatim}
\allowdisplaybreaks
\mathtoolsset{showonlyrefs}

\newcommand{\D}{\nabla}
\newcommand{\dd}{\partial}
\newcommand{\norm}[1]{\left\lVert#1\right\rVert}

\newtheorem{Thm}{Theorem}[section]
\newtheorem{Lem}[Thm]{Lemma}

\newtheorem{Def}[Thm]{Definition}
\newtheorem{Prop}[Thm]{Proposition}

\def\R{\mathbb R}
\def\H{\mathcal H}
\def\N{\mathbb N}
\def\e{\varepsilon}

\def\Xint#1{\mathchoice
   {\XXint\displaystyle\textstyle{#1}}
   {\XXint\textstyle\scriptstyle{#1}}
   {\XXint\scriptstyle\scriptscriptstyle{#1}}
   {\XXint\scriptscriptstyle\scriptscriptstyle{#1}}
   \!\int}
\def\XXint#1#2#3{{\setbox0=\hbox{$#1{#2#3}{\int}$}
     \vcenter{\hbox{$#2#3$}}\kern-.5\wd0}}

\def\dashint{\Xint-}
\newcommand{\prodscal}[2]{\langle{#1},{#2}\rangle}

\title{Ahlfors-regularity for minimizers of a multiphase optimal design problem}

\author{L. Esposito, L. Lamberti and G. Pisante}

\begin{document}
\maketitle

\begin{abstract}
 \noindent 
We establish an Alhfors-regularity result for minimizers of a multiphase optimal design problem. It is a variant of the classical variational problem which involves a finite number of chambers $\mathcal{E}(i)$ of prescribed volume that partition a given domain $\Omega\subset\R^n$. The cost functional associated with a configuration $\left(\{\mathcal{E}(i)\}_i,u\right)$ is made up of the perimeter of the partition interfaces and a Dirichlet energy term, which is discontinuous across the interfaces. We prove that the union of the optimal interfaces is $(n-1)$-Alhfors-regular via a penalization method and decay estimates of the energy.
\end{abstract}


	
\makeatletter

\makeatother

\section{Introduction}

The problem of partitioning an open domain into regions with minimal interface has deep roots in both classical geometry and modern variational analysis. Formally, the goal is to partition an open set $\Omega \subset \mathbb R^n$  into a finite collection of disjoint subsets
$ \{ \mathcal{E}(i) \}_{i=1}^N$ such that their union covers $ \Omega $ (up to a set of measure zero) and they minimize the total interfacial energy, basically interpreted as the $(n-1)$-dimensional Hausdorff measure of the common boundaries. A celebrated example in two dimensions is the Honeycomb Conjecture, resolved by T. C. Hales \cite{H}, which states that the regular hexagonal tiling minimizes the total perimeter among all partitions of the plane into regions of equal area. \\
\indent Such partition problems generalize classical isoperimetric inequalities and are intimately connected to the theory of minimal surfaces. They also arise naturally in various applications, including immiscible fluid separation (see \cite{Leo, W}), and image segmentation (see \cite{ABr,CT,TC}).\\
Mathematically, the problem often involves minimizing an energy functional of the form
\begin{equation}
\mathcal{P}(\{\mathcal{E}(i)\};\Omega) = \sum_{i<j} \mathcal{H}^{n-1}(\partial^* \mathcal{E}(i)\cap \partial^* \mathcal{E}(j)\cap\Omega),
\end{equation}
subject to volume constraints $|\mathcal{E}(i)| = m_i $. Here, $ \partial^* \mathcal{E}(i) $ denotes the reduced boundary of $ \mathcal{E}(i) $ in the sense of geometric measure theory, which captures the essential structure of the interface between phases.\\
\indent In this paper, we focus on functionals that depend not only on the interfacial energy of a partition, but also on a bulk energy term. To motivate this setting, we refer to the classical problem of liquid droplets subjected to an external electric field, where the equilibrium configuration is typically determined by the competition between interfacial and bulk energies. The interfacial energy, often modeled as proportional to the surface area of the droplet, reflects the action of surface tension and tends to favor compact shapes, such as spheres, that minimize the surface area for a given volume. \\
\indent However, the presence of an external electric field introduces a nonlocal bulk energy that accounts for the interaction of the electric field with the dielectric properties of the droplet and the surrounding medium. This contribution is typically expressed through the Dirichlet energy of the electrostatic potential. The balance between bulk energy and interfacial energy gives rise to a free boundary variational problem in which the domain itself is an unknown to be optimized. For a comprehensive study of this model, the reader can refer to the work of Muratov and Novaga, who have extensively analyzed the variational problems associated with charged liquid droplets, see \cite{DHV,MV,MN} and the references therein.\\
\indent A prototype version of functionals involving bulk and perimeter energies is the following:
\begin{equation}\label{model1}
\int_{\Omega}\sigma_{E}(x)|\nabla u|^2\,dx+ P(E;\Omega),
\end{equation}
with $u=u_0$ prescribed on $\partial \Omega$ and $\sigma_{E}(x)=\beta \mathbbm{1}_{E}+\alpha \mathbbm{1}_{\Omega\setminus E}$, $0<\alpha<\beta$.\\
This functional was formerly studied in {1993} in two papers by L. Ambrosio \& G. Buttazzo and F.H. Lin (see \cite{ABu, Lin}). Later on, refined regularity results for functionals of type \eqref{model1} have been established in \cite{DF,FJ} and for dimension two in \cite{LamLem,Lar1,Lar2}. Furthermore, the same problem has been studied in the case where both the bulk and interfacial energies are of a more general nature (see \cite{CEL1,CEL2,CFP,CFP2,E,EL,EL2,ELP,Lam,LK}).\\
\indent In this paper, we study optimal partitions associated with functionals that also depend on a bulk Dirichlet energy, which is discontinuous across the partition interface. To our knowledge, there are no regularity results in the literature for this context. The presence of multiple chambers significantly complicates the study of regularity due to the possibility of triple points or, even worse, multiple intersections between the chambers. \\
\indent Some notation is needed. Let $\Omega\subset \R^n$ a bounded connected open set and $N\in\N$ such that $N> 1$. An $N$-partition $\mathcal{E}$ of $\Omega$  is a family $\mathcal{E}=\{\mathcal{E}(i)\}_{i=1}^N$  of sets $\mathcal{E}(i)$ of finite perimeter with
\begin{equation*}
|\mathcal{E}(i)|>0,\quad\forall i\in\{1,\dots,N\},
\end{equation*}
\begin{equation*}
|\mathcal{E}(i)\cap\mathcal{E}(j)|=0,\quad\forall i,j\in\{1,\dots,N\},\, i<j.
\end{equation*}
\begin{equation*}
\sum_{i=1}^N|\mathcal{E}(i)|=|\Omega|.
\end{equation*}
We introduce the following main functional associated to a partition $\mathcal{E}$:
\begin{equation}\label{MainF}
{\mathcal F}(\mathcal{E},w)=\sum_{i=1}^N\int_{\mathcal{E}(i)}\alpha_i|\D w|^2\,dx +\frac{1}{2}\sum_{i=1}^NP(\mathcal{E}(i);\Omega),
\end{equation}
where the vector $\alpha=\left\{\alpha_i\right\}_{i=1}^{N}$ is positive, i.e. $\alpha_i>0$ for any $i\in\{1,\dots N\}$.\\
\indent  
 The interfaces of the $N$-partition $\mathcal{E}$ of $\Omega$ are the $\mathcal{H}^{n-1}$-rectifiable sets
$$\mathcal{E}(h, k) = \partial^* \mathcal{E}(h) \cap \partial^* \mathcal{E}(k)\cap \Omega,$$ 
where $0 \leq h, k \leq N$ and $h \neq k$.


Given $\{d_i\}_{i=1}^N$ such that 
$$d_i\in (0,|\Omega|),\, \forall i\in\{1,\dots,N\}\quad\text{and} \quad
\sum_{i=1}^N d_i=|\Omega|,
$$ 
we consider the minimization of the functional \eqref{MainF}
 assuming that the measures of the chambers $\mathcal{E}(i)$ are equal to $d_i$ and the function $w$ is prescribed on the boundary of $\Omega$. More precisely, given $u_0\in H^1(\Omega)$, we consider the following constrained problem:
\begin{equation}
\label{P_c}
\min
\left\{
\mathcal{F}(\mathcal{E},v)\,:\,  \mathcal{E}    \text{ is an } N\text{-partition of } \Omega,\,|\mathcal{E}_i|=d_i,\,i=1,\dots,N, \,v\in u_0+H_0^1(\Omega) 
\right\}.
\tag{$P_c$}
\end{equation}
The aim of the paper is to prove the $(n-1)$-Ahlfors-regularity of the interfaces of the optimal chambers. We recall that a  closed set $G\subset\R^n$ is said to be $(n-1)$-Ahlfors-regular if there exists a positive constant $C_A$ such that
\begin{equation*}
    C_A^{-1} r^{n-1}\leq\mathcal{H}^{n-1}(G\cap B_r(x_0))\leq C_A r^{n-1},\quad\forall x_0\in G,\,\forall r>0.
\end{equation*}
In particular, we prove the following theorem.

\begin{Thm}
\label{Main}
Let $(\mathcal{E},u)$ be a minimizer of the problem  $\eqref{P_c}$ and $U\Subset\Omega$ be an open set. Then, there exist a positive constant $C_A$ such that, for every $x_0\in\bigcup_{k=1}^N\dd\mathcal{E}(k)\cap\Omega$ and $B_r(x_0)\subset U$, it holds
\begin{equation*}
C_A^{-1} r^{n-1}\leq \sum_{k=1}^N P(\mathcal{E}(k);B_r(x_0))\leq C_A r^{n-1}.
\end{equation*}
Moreover, $\mathcal{H}^{n-1}\Big( {\Omega\cap}\bigcup_{k=1}^N\dd\mathcal{E}(k)\setminus \bigcup_{k=1}^N\dd^*\mathcal{E}(k)\Big)=0$ and $\bigcup_{k=1}^N\dd\mathcal{E}(k)$ is $(n-1)$-Ahlfors-regular.
\end{Thm}
The strategy of the proof of Theorem \ref{Main} follows a well-established path. First, we show that minimizers of \eqref{P_c} are indeed also minimizer for a penalized problem without constraint (see Theorem \ref{Teorema Penalizzazione}).
Afterwards, the proof follows by combining the upper and lower density estimates for the minimizers of the penalized problem contained in Theorem \ref{Upper Bound} and Theorem \ref{Lower Bound}.

\section{From constrained to penalized problem}
In the following theorem, we show that volume-constrained minimizers of $\eqref{P_c}$ are, in fact, unconstrained $\Lambda$-minimizers of the functional $\mathcal{F}$ defined in \eqref{MainF} (see Definition \ref{Lambda-Min} below). This type of relaxation of the volume constraint is standard in problems of this nature. To obtain this result, we employ a technique introduced in \cite{EF}, which, in our setting, is more intricate and requires a suitable adaptation due to the presence of multiple chambers.
\begin{Thm}
\label{Teorema Penalizzazione}
There exist $\Lambda_{0}>0$ such that  if $(\mathcal{E},u)$ is a minimizer of the functional
\begin{equation}
\label{Penalized}
{\mathcal F}_{\Lambda}(\mathcal{A},w)=\sum_{i=1}^N\int_{\mathcal{A}(i)}\alpha_i|\D w|^2\,dx +\frac{1}{2}\sum_{i=1}^NP(\mathcal{A}(i);\Omega)+\Lambda\sum_{i=1}^N\big||\mathcal{A}(i)| - d_i\big|,
\end{equation}
for some $\Lambda \geq \Lambda_0$, among all configurations $(\mathcal{A},w)$ such that $w=u_0$ on $\partial \Omega$,
then $|\mathcal{E}|=d$ and $(\mathcal{E},u)$ is a minimizer of problem \eqref{P_c}.
Conversely, if $(\mathcal{E},u)$ is a minimizer of problem \eqref{P_c} among all configurations $(\mathcal{A},w)$ such that $w=u_0$ on $\partial \Omega$, then it is a minimizer of \eqref{Penalized}, for all $\Lambda \geq \Lambda_0$.
\end{Thm}
\begin{proof}
The first part of the theorem can be proved by contradiction. We assume that there exist a positive sequence $(\Lambda_h)_{h\in \mathbb N}$ such that $\Lambda_h\rightarrow +\infty$, as $h\rightarrow +\infty$, and a sequence of configurations $(\mathcal{E}_h,u_h)$ minimizing $\mathcal{F}_{\Lambda_h}$ and such that $u_h=u_0$ on $\partial \Omega$ and $|\mathcal{E}_h|\neq d$, for all $h\in\N$. We choose an arbitrary fixed partition $\mathcal{E}_0$ of $\Omega$ such that $|\mathcal{E}_0|=d$. We point out that 
\begin{equation}\label{Theta}
\mathcal{F}_{\Lambda_h}(\mathcal{E}_h,u_h)\leq\mathcal{F}(\mathcal{E}_0,u_0):=\Theta.
\end{equation}
Our aim is to show that there exists a configuration $(\mathcal{\widetilde{E}}_h,\tilde{u}_h)$ such that, for $h$ sufficiently large, $\mathcal{F}_{\Lambda_h}(\mathcal{\widetilde{E}}_h,\tilde{u}_h)< \mathcal{F}_{\Lambda_h}({\mathcal{E}_h},{u_h})$, thus proving the result by contradiction.\\
\indent By condition $\eqref{Theta}$, it follows that the sequence $(u_h)_{h\in\N}$ is bounded in $H^1(\Omega)$, the perimeter of the partition $\mathcal{E}_h$ in $\Omega$ is uniformly bounded and $|\mathcal{E}_h(i)|\rightarrow d_i$, for any $i\in\{1,\dots,N\}$. Therefore, possibly extracting a not relabelled subsequence, we may assume that there exists a configuration $(\mathcal{E},u)$ such that $u_h$ converges to $u$ weakly in $H^1(\Omega)$, $\mathbbm{1}_{\mathcal{E}_h(i)}\rightarrow\mathbbm{1}_{\mathcal{E}(i)}$ a.e. in $\Omega$, where the collection $\mathcal{E}=\{\mathcal{E}(i)\}_{i=1}^N$ is a partition of $\Omega$ and $|\mathcal{E}|=d$. The couple $(\mathcal{E},u)$ will be used as a reference configuration for the definition of $(\widetilde{\mathcal{E}}_h,\tilde{u}_h)$.\\
\indent 
By appropriately rearranging the order of the chambers, we can assume that there exists an $i\in\{1,\dots,N\}$ such that $|\mathcal{E}_h(i)|<d_i$, for any $h\in\N$.
Since 
\begin{equation}
\sum_{j=1}^N|\mathcal{E}_h(j)|=|\Omega|=\sum_{j=1}^N d_j,
\end{equation}
we can also assume that there exists $j\in\{1,\dots,N\}$ such that $|\mathcal{E}_h(j)|>d_j$, for any $h\in\N$.\\
Let $0\leq i, j \leq N$, $i\neq j$. We say that $\mathcal{E}(i)$ and $\mathcal{E}(j)$ are
neighboring chambers, if $\mathcal{H}^{n-1}(\mathcal{E}(i, j)) > 0$. 
If there exist two neighboring chambers $\mathcal{E}(i)$ and $\mathcal{E}(j)$ with $|\mathcal{E}_h(i)|<d_i$ and $|\mathcal{E}_h(j)|>d_j$, then we can argue exactly as in \cite{EF} to construct the configuration $(\widetilde{\mathcal{E}}_h,\tilde{u}_h)$. Otherwise we will work with the pair of chambers $\mathcal{E}(\overline i)$ and $\mathcal{E}(\overline  j)$ with with $|\mathcal{E}_h(\overline  i)|<d_{\overline  i}$ and $|\mathcal{E}_h(\overline  j)|>d_{\overline  j}$ that are the \emph{closest} in a suitable sense. More precisely, for $i,j\in\{1,\dots,N\}$, we denote by $c_{ij}$ the order of link between the chambers $\mathcal{E}(i)$ and $\mathcal{E}(j)$ that is defined as the minimum number $m$ such that there exist chambers $\mathcal{E}(k_1),\cdots \mathcal{E}(k_m)$, such that $\mathcal{E}(i)$ is neighboring $\mathcal{E}(k_1)$, $\mathcal{E}(j)$ is neighboring $\mathcal{E}(k_m)$ and $\mathcal{E}(k_l)$ is neighboring $\mathcal{E}(k_{l+1})$ for any $l\in \{1,\dots,m-1\}$.
We identify $\overline{i}$ and $\overline{j}$ as the indices of two chambers such that 
\begin{equation}
(\overline{i},\overline{j})=\mathrm{argmin}\{c_{ij}\,:\,(i,j)\in\{1,\dots,N\},\,|\mathcal{E}_h(i)|<d_i,\,|\mathcal{E}_h(j)|>d_j\}.
\end{equation}
Therefore there exists $m\in\N$ such that $\mathcal{E}_h(i)$ and $\mathcal{E}_h(j)$ are linked through some chambers $\mathcal{E}_h(k_1),\dots,\mathcal{E}_h(k_m)$, where $k_\ell\in\{1,\dots,N\}\setminus\{\overline{i},\overline{j}\}$ and $|\mathcal{E}_h(k_\ell)|=d_{k_\ell}$, for any $\ell\in\{1,\dots,N\}$.
\\
We may also assume that there is only one intermediate chamber. i.e. $m=1$ or equivalently that $c_{\overline{ij}}=1$ as in the other case the construction can be carried over in a similar way. \\
\indent Then, for simplicity, up to relabeling the chambers, we will assume that $i=1$, $j=3$, $c_{13}=1$ and $\mathcal{E}_h(2)$ is the linking chamber between $\mathcal{E}_h(1)$ and $\mathcal{E}_h(3)$, with
\begin{equation*}
|\mathcal{E}_h(1)|<d_1,\quad |\mathcal{E}_h(2)|=d_2,\quad|\mathcal{E}_h(3)|>d_3.
\end{equation*}

\textbf{Step 1.} \emph{Construction of $(\tilde{\mathcal{E}}_h,\tilde{u}_h)$}. Let us choose $(\sigma_1)_h\in\R$ and $(\sigma_2)_h$ such that
\begin{equation}
\label{sigma1}
(\sigma_1)_h\in\bigg(0,\alpha\min\bigg\{\frac{1}{2^n},\frac{|d_1-|\mathcal{E}_h(1)||}{n2^{n+1}}\bigg\}\bigg),
\end{equation}
\begin{equation}
\label{sigma2}
(\sigma_2)_h\in\bigg(0,\min\bigg\{\frac{1}{2^n},\frac{|d_2-|\mathcal{E}_h(2)||}{n2^{n+1}}\bigg\}\bigg),
\end{equation}
where $\alpha=\alpha(n,N)\in(0,1)$ is a constant that will be chosen later. We fix $\ell\in\{1,2\}$. Since the three chambers are linked, the set $\dd^*\mathcal{E}(\ell)\cap\dd^*\mathcal{E}(\ell+1)\cap\Omega$ is not empty. Thus, we can take a point $x_{\ell}\in\dd^*\mathcal{E}(\ell)\cap\dd^*\mathcal{E}(\ell+1)\cap\Omega$ such that
\begin{equation*}
\lim_{\rho\rightarrow0^+}\frac{|\mathcal{E}(\ell)\cap B_\rho(x_{\ell})|}{\omega_n\rho^n}=\lim_{\rho\rightarrow0^+}\frac{|\mathcal{E}(\ell+1)\cap B_\rho(x_{\ell})|}{\omega_n\rho^n}=\frac{1}{2}.
\end{equation*}
Since
\begin{equation*}
1=\frac{|B_\rho(x_{\ell})|}{\omega_n\rho^n}=\frac{|\mathcal{E}(\ell)\cap B_\rho(x_{\ell})|}{\omega_n\rho^n}+\frac{|\mathcal{E}(\ell+1)\cap B_\rho(x_{\ell})|}{\omega_n\rho^n}+\sum_{\substack{k=1\\k\neq \ell,\ell+1}}^N\frac{|\mathcal{E}(k)\cap B_r(x_{\ell})|}{\omega_n\rho^n},
\end{equation*}
there exists $\eta\in(0,\min\big\{\frac{\mathrm{dist}(x_1,x_2)}{2},1\big\})$ such that if $0<r<\eta$, then
\begin{equation}
\label{q1}
\sum_{\substack{j=1\\j\neq \ell,\ell+1}}^N|\mathcal{E}(j)\cap B_r(x_{\ell})|<\kappa\min_{i=1,2}(\sigma_i)_h r^n,
\end{equation}
for some positive constant $\kappa$ that will be chosen later and $h$ sufficiently large.
By De Giorgi structure theorem for sets of finite perimeter, $\mathbbm{1}_{\frac{\mathcal{E}(\ell)-x_{\ell}}{r}}\rightarrow\mathbbm{1}_{H_{\ell}}$ in $L^1_{loc}(\R^n)$, as $r\rightarrow 0^+$, where $H_{\ell}=\{\prodscal{z}{\nu_{\mathcal{E}(\ell)}(x_l)}<0\}=\{\prodscal{z}{\nu_{\mathcal{E}(\ell+1)}(x_l)}>0\}$. Let $y_{\ell}\in B_1\setminus H_{\ell}$ be the point 
\begin{equation*}
y_{\ell}=\frac{\nu_{\mathcal{E}(\ell)}(x_l)}{2}=-\frac{\nu_{\mathcal{E}(\ell+1)}(x_l)}{2}.
\end{equation*}
Therefore, there exists $0<r<\eta$ such that \eqref{q1} holds and
\begin{equation}
\label{meas}
\left|\frac{\mathcal{E}(\ell)-x_{\ell}}{r}\cap B_{1/2}(y_{\ell})\right|<\e,\quad\left|\frac{\mathcal{E}(\ell)-x_{\ell}}{r}\cap B_1(y_{\ell})\right|>\frac{1}{2^{n+2}}.
\end{equation}
Setting $x_{\ell}':=x_{\ell}+ry_{\ell}$, from the convergence of $\mathcal{E}_h$ to $\mathcal{E}$, we have that, for $h$ sufficiently large,
\begin{equation}\label{unodue}
|\mathcal{E}_h(\ell)\cap B_{r/2}(x'_{\ell})|<\e r^n,\quad|\mathcal{E}_h(\ell)\cap B_{r}(x'_{\ell})|>\frac{ r^n}{2^{n+2}},\quad \sum_{\substack{j=1\\j\neq \ell,\ell+1}}^N|\mathcal{E}_h(j)\cap B_r(x_{\ell})|<\kappa\min_{i=1,2}(\sigma_i)_hr^n,
\end{equation}
where $\kappa$ will be chosen later. We remark that, since $r<\eta$, $B_r(x_1)$ and $B_r(x_2)$ are disjoint.\\
\indent Now we define the following bi-Lipschitz map used in \cite{EF} which maps $B_r(0)$ into itself: 
\begin{equation}\label{unotre}
\Psi(x):=
\begin{cases}
\bigl[1-(\sigma_\ell)_h\big(2^n-1\big)\bigr]x \quad&\text{if }|x|<\displaystyle\frac{r}{2},\\
x+(\sigma_\ell)_h\bigg(\displaystyle 1-\frac{r^n}{|x|^n}\bigg)x &\text{if }\displaystyle\frac{r}{2}\leq|x|<r,\\
x & \text{if }|x|\geq r.
\end{cases}
\end{equation}
We denote the corresponding action localized in the ball $B_r(x'_\ell)$ by
\begin{equation*}
\Phi_\ell(x)=x'_\ell+\Psi(x-x'_\ell),\quad\forall x\in\R^n.
\end{equation*}
In the sequel, we will use some estimates for the map $\Phi_{\ell}$ that can be easily obtained by direct computations (see \cite{EF} for the explicit calculations). These estimates are trivial for $|x-x'_\ell|<r/2$, whereas they
can be deduced by the explicit expression of $\nabla \Psi$ for $r/2<|x-x'_\ell|<r$, that is
\begin{equation*}
\frac{\partial (\Phi_{\ell})_i}{\partial x_j}(x)=\delta_{ij}+(\sigma_\ell)_h\bigg[\bigg(1-\frac{r^n}{|x-x'_\ell|^n}\bigg)\delta_{ij}+nr^n\frac{(x-x'_\ell)_i(x-x'_\ell)_j}{|x-x'_\ell|^{n+2}}\bigg],\quad \forall i,j\in\{1,\dots, n\}.
\end{equation*}
There exists a positive constant $C=C(n)$ depending only on $n$ such that, 
\begin{equation}\label{phi}
\norm{\nabla \Phi_{\ell}^{-1}(y)-I}\leq C(n)(\sigma_\ell)_h, \quad\forall y\in B_{r}(x'_{\ell}),
\end{equation}
\begin{equation}
\label{Jphi1}
1+C(n)(\sigma_\ell)_h\leq J\Phi_{\ell}(x)\leq1+n2^n (\sigma_\ell)_h,\quad\forall x\in B_r(x'_{\ell})\setminus B_{\frac{r}{2}}(x'_{\ell}).
\end{equation}
\begin{equation}
\label{Jphi2}
1-n(2^n-1)(\sigma_\ell)_h\leq J\Phi_{\ell}(x)\leq 1-(2^n-1)(\sigma_\ell)_h,\quad\forall x\in B_{\frac{r}{2}}(x'_{\ell}).
\end{equation}
Accordingly, for $j\neq\ell,\ell+1$, we estimate,
\begin{align}
\label{q5}
&\Big||\Phi_{\ell}(\mathcal{E}_h(j)\cap B_r(x'_{\ell}))|-|\mathcal{E}_h(j)\cap B_r(x'_{\ell})|\Big|\notag\\
& \leq\int_{\mathcal{E}_h(j)\cap B_r(x'_{\ell})}|J\Phi_{\ell}-1|\,dx\leq  n2^n(\sigma_\ell)_h|\mathcal{E}_h(j)\cap B_r(x'_{\ell})|\leq n2^n\kappa\min_{j\in\{1,2\}}(\sigma_j)_hr^n,
\end{align}
whenever, for $j=\ell,\ell+1$, 
\begin{align}
\label{q5b}
&\Big||\Phi_{\ell}(\mathcal{E}_h(j)\cap B_r(x'_{\ell}))|-|\mathcal{E}_h(j)\cap B_r(x'_{\ell})|\Big|\notag\\
& \leq\int_{\mathcal{E}_h(j)\cap B_r(x'_{\ell})}|J\Phi_{\ell}-1|\,dx\leq  n2^n(\sigma_\ell)_h|\mathcal{E}_h(j)\cap B_r(x'_{\ell})|\leq n2^n(\sigma_\ell)_h r^n.
\end{align}
We define
\begin{equation}\label{perturbation}
\Phi:=\Phi_{1}\circ\Phi_{2},\quad\tilde{\mathcal{E}}_h:=\{\Phi(\mathcal{E}_h(i))\}_{i=1}^N,\quad\tilde{u}_h:=u_h\circ\Phi^{-1},
\end{equation}
and remark that $\Phi$ act in not trivial way only in $B_r(x'_{1})$ and $B_r(x'_{2})$, leaving anything unchanged outside these balls.\\
\indent \textbf{Step 2. }\emph{$\Phi$ does not modify too much $|\mathcal{E}_h(1)|$ and $|\mathcal{E}_h(3)|$}. \\
Let us show that under conditions \eqref{sigma1} , \eqref{sigma2} on $(\sigma_1)_h$ and $(\sigma_2)_h$, it results
\begin{equation}
\label{q3}
|\Phi(\mathcal{E}_h(1))|<d_1,
\end{equation}
\begin{equation}
\label{q4}
|\Phi(\mathcal{E}_h(3))|>d_3.
\end{equation}
Since the application of the map $\Phi$ leaves the measure unchanged outside $B_r(x'_1)\cup B_r(x'_2)$, we can evaluate the differences in measure solely within $B_r(x'_1)\cup B_r(x'_2)$.
We have that
\begin{align}
& \Big||\Phi(\mathcal E_h(1))|-|\mathcal E_h(1)| \Big| \\
& \leq\Big||\Phi_1(\mathcal E_h(1)\cap B_r(x'_1))|-|\mathcal E_h(1)\cap B_r(x'_1)|\Big|+\Big||\Phi_2(\mathcal E_h(1)\cap B_r(x'_2))|-|\mathcal E_h(1)\cap B_r(x'_2)|\Big|.
\end{align}
Applying \eqref{q5} and \eqref{q5b}, by the choice of $(\sigma_1)_h$ we deduce
\begin{align*}
|\Phi(\mathcal E_h(1))|\leq |\mathcal E_h(1)|+n2^{n+1}(\sigma_1)_h<d_1.
\end{align*}
The same argument can be applied to chamber $\mathcal E_h(3)$ under the transformation by the map
$\Phi$; that is, we have
\begin{align}
&\Big||\Phi(\mathcal E_h(3))|-|\mathcal E_h(3)| \Big| \\
& \leq\Big||\Phi_1(\mathcal E_h(3)\cap B_r(x'_{1}))|-|\mathcal E_h(3)\cap B_r(x'_{1})|\Big|+\Big||\Phi_2(\mathcal E_h(3)\cap B_r(x'_{2}))|-|\mathcal E_h(3)\cap B_r(x'_{2})|\Big|.
\end{align}
Using again \eqref{q5} and \eqref{q5b}, by the choice of $(\sigma_2)_h$ we deduce
\begin{align*}
|\Phi(\mathcal E_h(3))|\geq |\mathcal E_h(3)|-n2^{n+1}(\sigma_2)_h>d_3.
\end{align*}
\textbf{Step 3. }\emph{The choice of $(\sigma_2)_h$}. In this step we prove that for every $(\sigma_1)_h$ as in \eqref{sigma1} there exists $(\sigma_2)_h$ as in \eqref{sigma2} such that 
\begin{equation}\label{q11}
|\Phi(\mathcal{E}_h(2))|=d_2.
\end{equation}
To this end, we prove that the action of the map $\Phi$ results in an increase in the measure of $\mathcal{E}_h(2)\cap B_r(x_{2}')$ and a decrease in the measure of $\mathcal{E}_h(2)\cap B_r(x_{1}')$, modulated by the parameters $(\sigma_2)_h$ and $(\sigma_1)_h$.
Precisely we prove that there exist two positive constants $C_g=C_g(n)$ and $C_l=C_l(n)$ such that
\begin{align}
  & \label{Gain}\mathcal{G}_h:=|\Phi_{2}(\mathcal{E}_h(2))\cap B_r(x_{2}')|-|\mathcal{E}_h(2)\cap B_r(x_{2}')|\geq C_g(n)(\sigma_2)_hr^n,\\
  & \label{loss}\mathcal{L}_h:=|\mathcal{E}_h(2)\cap B_r(x_{1}')|-|\Phi_{1}(\mathcal{E}_h(2))\cap B_r(x_{1}')|\geq C_l(n)(\sigma_1)_hr^n.
\end{align}
With this notations \eqref{q11} can be rephrased
\begin{equation*}
    \mathcal{G}_h=\mathcal{L}_h
\end{equation*}
Inequality \eqref{Gain} is the simplest to prove. Indeed, taking the estimates \eqref{Jphi1} and \eqref{Jphi2}  on $J\Phi_{2}$ and \eqref{unodue} into account, we have that
\begin{align}\label{gain}
& \mathcal{G}_h=|\Phi_{2}(\mathcal{E}_h(2))\cap B_{r}(x'_{2})|-|\mathcal{E}_h(2)\cap B_r(x_{2}')|=\int_{\mathcal{E}_h(2)\cap B_r(x_{2}')}(J\Phi_{2}-1)\,dx\\
&=\int_{\mathcal{E}_h(2)\cap B_{\frac{r}{2}}(x_{2}')}(J\Phi_{2}-1)\,dx+
\int_{\mathcal{E}_h(2)\cap (B_{r}(x_{2}')\setminus B_{\frac{r}{2}}(x_{2}'))}(J\Phi_{2}-1)\,dx\\
&\geq -\varepsilon n\big(2^n-1\big)(\sigma_2)_hr^n+C(n)\bigg(\frac{\omega_n}{2^{n+2}}-\varepsilon\bigg)(\sigma_2)_hr^n\\
& =\bigg(\big[-n\big(2^n-1\big)-C(n)\big]\varepsilon+C(n)\frac{\omega_n}{2^{n+2}}\bigg)(\sigma_2)_hr^n={C}_g(n)(\sigma_2)_hr^n,
\end{align}
where $\varepsilon$ is chosen small enough in such a way that ${C}_g(n)>0$.\\
\indent To prove \eqref{loss}, we first note that by employing exactly the same computations used to establish \eqref{Gain}, we deduce that $\mathcal{E}_h(1)$ increases in measure within $B_r(x_{1}')$. Specifically, we have that
\begin{equation}\label{Gain1}
 |\Phi_{1}(\mathcal{E}_h(1))\cap B_r(x_{1}')|-|\mathcal{E}_h(1)\cap B_r(x_{1}')|\geq {C_g}(n)(\sigma_1)_hr^n.
\end{equation}
Since the total measure of the ball $B_r(x_{1}')$ is preserved by the map $\Phi_1$ we have that
\begin{equation}
 \sum_{j=1}^N  |\Phi_{1}(\mathcal{E}_h(j))\cap B_r(x_{1}')|  =\sum_{j=1}^N|\mathcal{E}_h(j))\cap B_r(x_{1}')|,
\end{equation}
consequently we deduce
\begin{align}
&|\mathcal{E}_h(2)\cap B_r(x_{1}')|-|\Phi_{1}(\mathcal{E}_h(2))\cap B_r(x_{1}')|\\
&=|\Phi_{1}(\mathcal{E}_h(1))\cap B_r(x_{1}')|-|\mathcal{E}_h(1)\cap B_r(x_{2}')|+\sum_{\substack{j=1\\j\neq 1,2}}^N\big[|\Phi_{1}(\mathcal{E}_h(j))\cap B_r(x_{1}')|-|\mathcal{E}_h(j)\cap B_r(x_{1}')|\big]\\
&\geq {C_g}(n)(\sigma_1)_hr^n -\kappa n(N-2)2^n(\sigma_1)_hr^n=[{C_g}(n)-\kappa n(N-2)2^n](\sigma_1)_hr^n,
\end{align}
where we used \eqref{Gain1} and \eqref{q5}.
Therefore \eqref{loss} is proved with $C_l(n)={C_g}(n)/2$ if we chose
$\kappa={C_g}(n) / \big[n(N-2) 2^{n+1}\big]$. If we denote with $C_g^{m}$ and $C_l^{m}$ the greatest constant such that \eqref{Gain} and \eqref{loss} holds true we deduce that
\begin{align}
  & \label{Gainf}|\Phi_{2}(\mathcal{E}_h(2))\cap B_r(x_{2}')|-|\mathcal{E}_h(2)\cap B_r(x_{2}')|= C_g^{m}(n)(\sigma_2)_hr^n\\
  & \label{lossf}|\mathcal{E}_h(2)\cap B_r(x_{1}')|-|\Phi_{1}(\mathcal{E}_h(2))\cap B_r(x_{1}')|= C_l^{m}(n)(\sigma_1)_hr^n.
\end{align}
Finally we can conclude observing that
\begin{equation*}
\mathcal{G}_h-\mathcal{L}_h=[C_g^{m}(n)(\sigma_2)_h-C_l^{m}(n)(\sigma_1)_h]r^n.
\end{equation*}
Then \eqref{q11} is proved if we choose $(\sigma_2)_h=\frac{C_l^{m}(n)}{C_g^{m}(n)}(\sigma_1)_h$.

\textbf{Step 4: }\emph{Reaching a contradiction.} 
Finally, we prove that the perturbation defined in \eqref{perturbation} leads to a decrease in energy, namely  
\begin{equation*}
{\mathcal F}_{\Lambda_h}(\mathcal{E}_h,u_h)-{\mathcal F}_{\Lambda_h}(\tilde{\mathcal{E}}_h,\tilde{u}_h) > 0,
\end{equation*}
which contradicts the minimality of $(\mathcal{E}_h, u_h)$.
For convenience of notation, we reformulate the energy difference by setting
\begin{align}\label{unoquattro}
& {\mathcal F}_{\Lambda_h}(\mathcal{E}_h,u_h)-{\mathcal F}_{\Lambda_h}(\tilde{\mathcal{E}}_h,\tilde{u}_h) \notag\\
& = \sum_{i=1}^N\sum_{\ell=1}^2\biggl[\int_{B_{r}(x'_{\ell})}\sigma_{\mathcal{E}_h(i)}|\D u_h|^2\,dx-\int_{B_{r}(x'_{\ell})}\sigma_{\Phi_{\ell}(\mathcal{E}_h(i))}|\D( u_h\circ\Phi_{\ell}^{-1})|^2\,dx\bigg]\notag\\
& +\frac{1}{2}\sum_{i=1}^N\sum_{\ell=1}^{2}\bigl[P(\mathcal{E}_h(i),{\overline B}_{r}(x'_{\ell}))-P(\Phi_{\ell}(\mathcal{E}_h(i)),{\overline B}_{r}(x'_{\ell}))\bigr]\notag\\
& +\Lambda_h\sum_{i=1}^N[||\mathcal{E}_h(i)|-d_i|-||\Phi(\mathcal{E}_h(i))|-d_i|]= I_{1,h}+I_{2,h}+I_{3,h}. 
\end{align}
{\bf Substep 1.} {\em Estimate of $I_{1,h}$}. Performing the change of variables $y=\Phi_{\ell}(x)$, and observing that $\mathbbm{1}_{\Phi(\mathcal{E}_h(i))}\circ \Phi=\mathbbm{1}_{\mathcal{E}_h(i)}$, we get
\begin{equation}
I_{1,h}=\sum_{i=1}^N\sum_{\ell=1}^2\int_{B_r(x'_{\ell})}\sigma_{\mathcal{E}_h(i)}(x)\bigl[|\nabla u_h(x)|^2-\bigl|\nabla u_h(x)\circ\nabla \Phi_{\ell}^{-1}(\Phi_{\ell}(x))\bigl|^2J\Phi_{\ell}(x)\bigr]\,dx=A_{1,h}+A_{2,h},
\end{equation}
where $A_{1,h}$ stands for the above integral over $B_{\frac{r}{2}}(x'_{\ell})$ and $A_{2,h}$ for the same integral over $B_{r}(x'_{\ell})\setminus B_{\frac{r}{2}}(x'_{\ell})$. Using the explicit expressions of $\nabla \Phi_{\ell}^{-1}$ and $J\Phi$,  which remain constant inside $B_{\frac{r}{2}}(x'_{\ell})$, we easily get
\begin{align*}
A_{1,h}=\sum_{i=1}^N\sum_{\ell=1}^2\int_{B_{\frac{r}{2}}(x'_{\ell})}\sigma_{\mathcal{E}_h(i)}|\D u_h|^2\Big\{1-\big[1-(\sigma_\ell)_h(2^n-1)\big]^{n-2}\Big\}\,dx\geq 0.
\end{align*}
Inside $B_{r}(x'_{\ell})\setminus B_{\frac{r}{2}}(x'_{\ell})$, even though  $\nabla\Phi_{\ell}^{-1}$ is not constant, we can use  \eqref{phi} and \eqref{Jphi1} to obtain
\begin{align}
A_{1,h}
& \geq \sum_{i=1}^N\sum_{\ell=1}^2\int_{B_{r}(x'_{\ell})\setminus B_{\frac{r}{2}}(x'_{\ell})}\sigma_{\mathcal{E}_h(i)}|\D u_h|^2\Big\{1-\big[1-	(\sigma_\ell)_h\big(2^n-1\big)\big]^{-2}\big(1+2^nn(\sigma_\ell)_h\big)\Big\}\,dx\\
& \geq -c(n)\sum_{i=1}^N\sum_{\ell=1}^2(\sigma_\ell)_h\int_{B_{r}(x'_{\ell})\setminus B_{\frac{r}{2}}(x'_{\ell})}\sigma_{\mathcal{E}_h(i)}|\D u_h|^2\,dx\\
& \geq-c(n)\Theta((\sigma_1)_h+(\sigma_2)_h),
\end{align}
where we also used \eqref{Theta},
thus getting
\begin{equation}
\label{I1h}
I_{1,h}\geq -\overline{C}_1\Theta((\sigma_1)_h+(\sigma_2)_h),
\end{equation}
for some positive constant $\overline{C}_1=\overline{C}_1(n)$.\\
\indent {\bf Substep 2.} {\em Estimate of $I_{2,h}$}. In order to estimate $I_{2,h}$, we can use the area formula for maps between rectifiable sets. If we denote by $T^{\ell,i}_{h,x}$ the tangential gradient of $\Phi_{\ell}$ along the approximate tangent space to $\partial^* \mathcal{E}_h(i)$ in $x$ and $\big(T^{\ell,i}_{h,x}\big)^*$ is the adjoint of the map $T^{\ell,i}_{h,x}$, the $(n-1)$-dimensional jacobian of $T^{\ell,i}_{h,x}$ is given by
\begin{equation}
J_{n-1}T^{\ell,i}_{h,x}=\sqrt{{\rm det}\bigl(\big(T^{\ell,i}_{h,x}\big)^*\circ T^{\ell,i}_{h,x}\bigr)}.
\end{equation}
Thereafter we can estimate
\begin{equation}\label{TJ}
J_{n-1}T^{\ell,i}_{h,x}\leq 1+(\sigma_\ell)_h+2^n(n-1)(\sigma_\ell)_h.
\end{equation}
We address the reader to \cite{EF} where explicit calculations are given. In order to estimate $I_{2,h}$, we use the area formula for maps between rectifiable sets (\cite[Theorem~2.91]{AFP}), thus getting
\begin{align*}
I_{2,h}
& =\sum_{i=1}^N\sum_{\ell=1}^2\big[P(\mathcal{E}_h(i);{\overline B}_r(x'_{\ell}))-P(\Phi_{\ell}(\mathcal{E}_h);{\overline B}_r(x'_{\ell}))\big]\\
& = \sum_{i=1}^N\sum_{\ell=1}^2\bigg[\int_{\partial^*\mathcal{E}_h(i)\cap{\overline B}_r(x'_{\ell})}d\H^{n-1}-\int_{\partial^*\mathcal{E}_h(i)\cap{\overline B}_r(x'_{\ell})}J_{n-1}T^{\ell,i}_{h,x}\,d\H^{n-1}\bigg] \\
& =\sum_{i=1}^N\sum_{\ell=1}^2\bigg[\int_{\partial^*\mathcal{E}_h(i)\cap{\overline B}_r(x'_{\ell})\setminus B_{\frac{r}{2}}(x'_{\ell})}\left(1-J_{n-1}T^{\ell,i}_{h,x}\right)\,d\H^{n-1}\\
& +\int_{\partial^*\mathcal{E}_h(i)\cap{\overline B}_{\frac{r}{2}}(x'_{\ell})}\left(1-J_{n-1}T^{\ell,i}_{h,x}\right)\,d\H^{n-1}	\bigg].
\end{align*}
Notice that the last integral in the above formula is non-negative since $\Phi_{\ell}$ is a contraction in $B_{\frac{r}{2}}$, hence $J_{n-1}T^{\ell,i}_{h,x}<1$ in $B_{r/2}$, while from \eqref{TJ} we have
\begin{align}
& \sum_{i=1}^N\sum_{\ell=1}^2\int_{\partial^*\mathcal{E}_h(i)\cap{\overline B}_r(x'_{\ell})\setminus B_{\frac{r}{2}}(x'_{\ell})}\left(1-J_{n-1}T^{\ell,i}_{h,x}\right)\,d\H^{n-1}\\
& \geq-2^nn((\sigma_1)_h+(\sigma_2)_h)\sum_{i=1}^N P(\mathcal{E}_h(i);{\overline B}_r(x'_{\ell}))\geq-2^nn\Theta((\sigma_1)_h+(\sigma_2)_h)\\
& =:-\overline{C}_2(n,\Theta)((\sigma_1)_h+(\sigma_2)_h),
\end{align}
thus concluding that
\begin{equation}\label{I2h}
I_{2,h}\geq-\overline{C}_2(n,\Theta)((\sigma_1)_h+(\sigma_2)_h).
\end{equation}
{\bf Substep 3.} {\em Estimate of $I_{3,h}$}. Taking  \eqref{q3}, \eqref{q4} and \eqref{q11} into account, we split $I_{3,h}$ in three addends as follows:
\begin{align}
\Lambda_h^{-1}\label{q7}
I_{3,h}=
& \big[|\Phi(\mathcal{E}_h(1))|-\mathcal{E}_h(1)|\big]+\big[|\mathcal{E}_h(3)|-|\Phi(\mathcal{E}_h(3))|\big]\notag\\
& +\sum_{i=4}^N\big[\big||\mathcal{E}_h(i)|-d_i\big|-\big||\Phi(\mathcal{E}_h(i))|-d_i\big|\big]=B_{1,h}+B_{2,h}+B_{3,h}.
\end{align}
We estimate the addends separately.  For $i>3$, by \eqref{q5} we have
\begin{align}
\label{q8}
& \big||\mathcal{E}_h(i)|-d_i|-\big||\Phi(\mathcal{E}_h(i))|-d_i\big|\leq \big||\mathcal{E}_h(i)|-|\Phi(\mathcal{E}_h(i))|\big|\notag\\
& \leq \big| |\mathcal{E}_h(i)\cap B_r(x'_{1})|-|\Phi_{1}(\mathcal{E}_h(i)\cap B_r(x'_{1}))|\big|+\big||\mathcal{E}_h(i)\cap B_r(x'_{2})|-|\Phi_{2}(\mathcal{E}_h(i)\cap B_r(x'_{2}))| \big|\notag\\
& \leq n2^{n+1}\kappa\min_{j\in\{1,2\}}(\sigma_j)_h r^n.
\end{align}
Therefore, we have
\begin{equation}\label{B3h}
B_{3,h}\geq -n(N-3) 2^{n+1}\kappa\min_{j\in\{1,2\}}(\sigma_j)_h r^n .  
\end{equation}
Regarding $B_{1,h}$, we use \eqref{Gain1} and \eqref{q5} to get
\begin{align}
\label{q9}
& B_{1,h}=|\Phi(\mathcal{E}_h(1))|-|\mathcal{E}_h(1)|\notag\\
& = \big[|\Phi_{1}(\mathcal{E}_h(1))\cap B_r(x'_{1})|-|\mathcal{E}_h(1)\cap B_r(x'_{1})|\big]+\big[|\Phi_{2}(\mathcal{E}_h(1))\cap B_r(x'_{2})|-|\mathcal{E}_h(1)\cap B_r(x'_{2})|\big]\notag\\
& \geq C_g(n)(\sigma_1)_h r^n-n2^n(\sigma_2)_h\kappa r^n.
\end{align}
Similarly, we can estimate $B_{2,h}$ in $B_r(x'_{2})$ observing that
the total measure of the ball $B_r(x_{2}')$ is preserved by the map $\Phi_2$, that is
\begin{equation}
 \sum_{j=1}^N  |\Phi_{2}(\mathcal{E}_h(j))\cap B_r(x_{2}')|  =\sum_{j=1}^N|\mathcal{E}_h(j))\cap B_r(x_{2}')|.
\end{equation}
Accordingly, using \eqref{gain} and \eqref{q5}, we deduce that
\begin{align}
& |\mathcal{E}_h(3)\cap B_r(x'_{2})|-|\Phi_{2}(\mathcal{E}_h(3))\cap B_r(x'_{2})|\\   
&=\big[|\Phi_{2}(\mathcal{E}_h(2)\cap B_r(x'_{2}))|-|\mathcal{E}_h(2)\cap B_r(x'_{2})|\big]+\sum_{j\neq2,3}\big[|\Phi_{2}(\mathcal{E}_h(j))\cap B_r(x'_{2})|-\mathcal{E}_h(j)\cap B_r(x'_{2})|\big]\\
&\geq C_g(n)(\sigma_2)_hr^n-n2^n(N-2)\kappa\min_{i\in\{1,2\}}(\sigma_i)_hr^n.
\end{align}
Therefore, we can conclude, using again \eqref{q5},
\begin{align}
\label{q10}
&B_{2,h}= |\mathcal{E}_h(3)|-|\Phi(\mathcal{E}_h(3))|\notag\\
& = \big[|\mathcal{E}_h(3)\cap B_r(x'_{2})|-|\Phi_{2}(\mathcal{E}_h(3))\cap B_r(x'_{2})|\big]+\big[|\mathcal{E}_h(3)\cap B_r(x'_{1})|-|\Phi_{1}(\mathcal{E}_h(3))\cap B_r(x'_{1})|\big]\notag\\
& \geq C_g(n)(\sigma_2)_hr^n-n2^n(N-2)\kappa\min_{i\in\{1,2\}}(\sigma_i)_hr^n-n2^n\kappa(\sigma_1)_hr^n\notag\\
& =C_g(n)(\sigma_2)_hr^n-n2^n(N-1)\kappa\min_{i\in\{1,2\}}(\sigma_i)_hr^n.
\end{align}
Finally, combining \eqref{q7}, \eqref{B3h}, \eqref{q9} and \eqref{q10}, we conclude that
\begin{align}
\label{I3h}
& I_{3,h}\geq \Lambda_h\big[C_g(n)-n2^n(N+1)\kappa\big]((\sigma_1)_h+(\sigma_2)_h)r^n.
\end{align}
\\
{\bf Substep 4.} {\em The contradiction}.
\noindent Inserting \eqref{I1h}, \eqref{I2h}, \eqref{I3h} in \eqref{unoquattro}, we conclude that
\begin{align*}
& I_{1,h}+I_{2,h}+I_{3,h} \geq [(\sigma_1)_h+(\sigma_2)_h]\big[-\overline{C}_1\Theta-\overline{C}_2+\Lambda_h\big(C_g(n)r^n-n2^n(N+1)\kappa r^n\big)\big]>0,
\end{align*}
if $\kappa=\kappa(n,N)$ is sufficiently small and $\Lambda_h\geq \Lambda_0=\Lambda_0(n,N,\Theta)$. This contradicts the minimality of $(\mathcal{E}_h,u_h)$, thus concluding the proof.
\end{proof}
The previous theorem  motivates the following definition.
\begin{Def}
[$\Lambda$-minimizers]\label{Lambda-Min}
The energy pair $(\mathcal{E},u)$ is a $\Lambda$-minimizer in $\Omega$ of the functional $\mathcal {F}$, defined in \eqref{MainF}, if and only if for every $B_r(x_0)\subset \Omega$ it holds that
\begin{equation}
    \mathcal{F}(\mathcal{E},u;B_r(x_0))\leq\mathcal{F}(\mathcal{G},v;B_r(x_0))+\Lambda\sum_{i=1}^N |\mathcal{G}(k)\Delta \mathcal{E}(k)|,
\end{equation}
whenever $(\mathcal{G},v)$ is an admissible test pair, namely, $\mathcal{G}$ is an $N$-partion of $\Omega$ such that $\mathcal{G}(k)\Delta \mathcal{E}(k)\subset \subset B_r(x_0)$ and $v-u\in H^1_0(B_r(x_0))$.
\end{Def}
\section{Energy decay estimates}
We start by proving a fundamental lemma which establishes that, for any ball that is either almost entirely contained within a single chamber or lies at the interface of only two chambers, the Dirichlet part of the functional satisfies a favorable decay estimate.
\begin{Lem}
\label{Lemma decadimento 1}
Let $(\mathcal{E},u)$ be a $\Lambda$-minimizer of the functional ${\mathcal F}$ defined in \eqref{MainF}. There exists $\tau_0\in(0,1)$ such that the following statement is true: for all $\tau \in (0,\tau_0)$ there exists $\varepsilon_1=\varepsilon_1(\tau)>0$ such that if $B_r(x_0)\subset \subset \Omega$ and one of the following conditions holds:
\begin{itemize}
\item[\emph{(i)}] There exists $i\in\{1,\dots,N\}$ such that $\frac{|\mathcal{E}(k)\cap B_r(x_0)|}{ |B_r(x_0)|}<\varepsilon_1$, for any $k\neq i$,
\item[\emph{(ii)}] There exist $i,j\in\{1,\dots,N\}$ and a half-space $H$ such that $\frac{\left|(\mathcal{E}(i)\setminus H)\cap B_r(x_0)\right|}{|B_r(x_0)|}<\varepsilon_1$,\break $\frac{\left|(\mathcal{E}(j)\cap H)\cap B_r(x_0)\right|}{|B_r(x_0)|}<\varepsilon_1$ and $\frac{|\mathcal{E}(k)\cap B_r(x_0)|}{|B_r(x_0)|}<\varepsilon_1$, for any $k\neq i,j$,
\end{itemize}
then
\begin{equation}\label{decayD}
\int_{B_{\tau r}(x_0)}|\nabla u|^2\,dx\leq C_1\tau ^{n} \int_{B_r(x_0)}|\nabla u|^2\,dx,
\end{equation}
for some positive constant $C_1=C_1(n,N,\alpha_k)$.
\end{Lem}

\begin{proof}
Let us fix $B_r(x_0)\subset \subset \Omega$ and $0<\tau<1$. Without loss of generality, we may assume that $\tau < 1/2$ and $x_0=0$. We start proving (i). We denote by $v$ the harmonic function in $B_{\frac{r}{2}}$ satisfying the condition $v=u$ on $B_{r}\setminus B_{\frac{r}{2}}$. Let $\phi\in H^1_0(B_{\frac{r}{2}})$. It holds that
\begin{equation}
\label{q12}
\int_{B_{\frac{r}{2}}}\prodscal{\D v}{\D\phi}\,dx=0.
\end{equation}
On the other hand, $u$ solves the following equation:
\begin{equation}
\label{q13}
\alpha_i\int_{B_{\frac{r}{2}}\cap\mathcal{E}(i)}\prodscal{\D u}{\D\phi}\,dx=-\sum_{k\neq i}\alpha_k\int_{B_{\frac{r}{2}}\cap\mathcal{E}(k)}\prodscal{\D u}{\D\phi}\,dx.
\end{equation}
Adding the equation
\begin{equation}
\alpha_i\int_{B_{\frac{r}{2}}\setminus\mathcal{E}(i)}\prodscal{\D u}{\D\phi}\,dx= \alpha_i\sum_{k\neq i}\int_{B_{\frac{r}{2}}\cap\mathcal{E}(k)}\prodscal{\D u}{\D\phi}\,dx,
\end{equation}
to \eqref{q13} and dividing by $\alpha_i$, we obtain
\begin{equation*}
\int_{B_{\frac{r}{2}}}\prodscal{\D u}{\D\phi}\,dx=\frac{1}{\alpha_i}\sum_{k\neq i}(\alpha_i-\alpha_k)\int_{B_{\frac{r}{2}}\cap\mathcal{E}(k)}\prodscal{\D u}{\D\phi}\,dx.
\end{equation*}
Now we choose $\phi:=v-u$ and we subtract \eqref{q12} from the previous equation, getting
\begin{align*}
\int_{B_{\frac{r}{2}}}|\D(v-u)|^2\,dx
& =\frac{1}{\alpha_i}\sum_{k\neq i}(\alpha_i-\alpha_k)\int_{B_{\frac{r}{2}}\cap\mathcal{E}(k)}\prodscal{\D u}{\D\phi}\,dx\\
& \leq \frac{\alpha}{\alpha_i} \sum_{k\neq i}\int_{B_{\frac{r}{2}}\cap \mathcal{E}(k)}|\prodscal{\D u}{\D\phi}|\,dx\\
& \leq \frac{\alpha}{\alpha_i}\Bigg(\sum_{k\neq i}\int_{B_{\frac{r}{2}}\cap \mathcal{E}(k)}|\D u|^2\,dx\Bigg)^{\frac{1}{2}}\bigg(\int_{B_{\frac{r}{2}}}|\D(v-u)|^2\,dx \bigg)^{\frac{1}{2}},
\end{align*}
where we have denoted $\alpha:=\max_{k}|\alpha_i-\alpha_k|$. 
Thus we infer that
\begin{equation*}
\int_{B_{\tau r}}|\D(v-u)|^2\,dx\leq \bigg(\frac{\alpha}{\alpha_i}\bigg)^2\sum_{k\neq i}\int_{B_{\frac{r}{2}}\cap \mathcal{E}(k)}|\D u|^2\,dx.
\end{equation*}
By the higher integrability for quadratic functional, (see for example \cite[Lemma 2.2]{FJ}, where an explicit calculation of the constants is provided)
\begin{equation*}
\dashint_{B_{\frac r2}}|\nabla u|^{2p}dx\leq C\Bigl(\dashint _{B_{r}}|\nabla u|^2 dx\Bigr)^p,
\end{equation*}
for some $p$ and $C$ both depending on $n,\alpha$.
Then using H\"older inequality, we deduce
\begin{align}
\int_{B_{\tau r}}|\D(v-u)|^2\,dx
& \leq C^{\frac{1}{p}}\bigg(\frac{\alpha}{\alpha_i}\bigg)^2\sum_{k\neq i}\bigg(\frac{|\mathcal{E}(k)\cap B_r|}{|B_r|}\bigg)^{1-\frac{1}{p}}\int_{B_r}|\D u|^2\,dx.
\end{align}
We choose $\varepsilon_1$ such that $\varepsilon_1^{1-\frac{1}{p}}=\tau^n$. Since $v$ is harmonic, we finally get
\begin{align*}
\int_{B_{\tau r}}|\D u|^2\,dx
& \leq 2\int_{B_{\tau r}}|\D (v-u)|^2\,dx+2\int_{B_{\tau r}}|\D v|^2\,dx\\
& \leq C^{\frac{1}{p}}\bigg(\frac{\alpha}{\alpha_i}\bigg)^2\tau^n\int_{B_r}|\D u|^2\,dx+2^{n+1}\tau^n \int_{B_r}|\D v|^2\,dx\\
& \leq C(n,N,\alpha,\alpha_i)\tau^n\int_{B_r}|\D u|^2\,dx,
\end{align*}
where we have also used the minimality of $v$.\\ 
\indent We are left with case (ii). We denote by $v$ the minimizer of the energy
\begin{equation*}
\int_{B_{\frac{r}{2}}}(\alpha_i\mathbbm{1}_{H}+\alpha_j\mathbbm{1}_{B_{\frac{r}{2}}\setminus H})|\D v|^2\,dx,
\end{equation*}
with the condition $v=u$ on $B_{r} \setminus B_{\frac{r}{2}}$. Let $\phi\in H^1_0(B_{\frac{r}{2}})$. It holds that
\begin{equation}
\label{q17}
\alpha_i\int_{B_{\frac{r}{2}}\cap H}\prodscal{\D v}{\D \phi}\,dx+\alpha_j\int_{B_{\frac{r}{2}}\setminus H}\prodscal{\D v}{\D \phi}\,dx=0.
\end{equation}
Now we rewrite the equation \eqref{q13}.
\begin{equation}
\label{q14}
\alpha_i\int_{B_{\frac{r}{2}}\cap \mathcal{E}(i)}\prodscal{\D u}{\D \phi}\,dx+\alpha_j\int_{B_{\frac{r}{2}}\cap \mathcal{E}(j)}\prodscal{\D u}{\D \phi}\,dx+\sum_{k\neq i,j}\alpha_k\int_{B_{\frac{r}{2}}\cap \mathcal{E}(k)}\prodscal{\D u}{\D \phi}\,dx=0.
\end{equation}
We decompose
\begin{align}
\label{q15}
& \alpha_i\int_{B_{\frac{r}{2}}\cap \mathcal{E}(i)}\prodscal{\D u}{\D \phi}\,dx=\alpha_i\int_{B_{\frac{r}{2}}\cap \mathcal{E}(i)\cap H}\prodscal{\D u}{\D \phi}\,dx+\alpha_i\int_{B_{\frac{r}{2}}\cap \mathcal{E}(i)\setminus H}\prodscal{\D u}{\D \phi}\,dx\notag\\
& =\alpha_i\int_{B_{\frac{r}{2}}\cap H}\prodscal{\D u}{\D \phi}\,dx+\alpha_i\bigg(\int_{B_{\frac{r}{2}}\cap \mathcal{E}(i)\setminus H}\prodscal{\D u}{\D \phi}\,dx\notag\\
& -\int_{B_{\frac{r}{2}}\cap\mathcal{E}(j)\cap H}\prodscal{\D u}{\D \phi}\,dx-\sum_{k\neq i,j}\int_{B_{\frac{r}{2}}\cap\mathcal{E}(k)\cap H}\prodscal{\D u}{\D \phi}\,dx\bigg)
\end{align}
and similarly
\begin{align}
\label{q16}
& \alpha_j\int_{B_{\frac{r}{2}}\cap \mathcal{E}(j)}\prodscal{\D u}{\D \phi}\,dx=\alpha_j\int_{B_{\frac{r}{2}}\setminus H}\prodscal{\D u}{\D \phi}\,dx+\alpha_j\bigg(\int_{B_{\frac{r}{2}}\cap \mathcal{E}(j)\cap H}\prodscal{\D u}{\D \phi}\,dx\notag\\
& -\int_{B_{\frac{r}{2}}\cap\mathcal{E}(i)\setminus H}\prodscal{\D u}{\D \phi}\,dx-\sum_{k\neq i,j}\int_{B_{\frac{r}{2}}\cap\mathcal{E}(k)\setminus H}\prodscal{\D u}{\D \phi}\,dx\bigg).
\end{align}
Inserting \eqref{q15} and \eqref{q16} in \eqref{q14}, we deduce that
\begin{align}
& \alpha_i\int_{B_{\frac{r}{2}}\cap H}\prodscal{\D u}{\D \phi}\,dx+\alpha_j\int_{B_{\frac{r}{2}}\setminus H}\prodscal{\D u}{\D \phi}\,dx\\
& =(\alpha_i-\alpha_j)\int_{B_{\frac{r}{2}}\cap \mathcal{E}(j)\cap H}\prodscal{\D u}{\D \phi}\,dx+(\alpha_j-\alpha_i)\int_{B_{\frac{r}{2}}\cap \mathcal{E}(i)\setminus H}\prodscal{\D u}{\D \phi}\,dx\\
& +\sum_{k\neq i,j}\bigg(\alpha_i\int_{B_{\frac{r}{2}}\cap \mathcal{E}(k)\cap H}\prodscal{\D u}{\D \phi}\,dx+\alpha_j\int_{B_{\frac{r}{2}}\cap \mathcal{E}(k)\setminus H}\prodscal{\D u}{\D \phi}\,dx-\alpha_k\int_{B_{\frac{r}{2}}\cap \mathcal{E}(k)}\prodscal{\D u}{\D \phi}\,dx\bigg).
\end{align}
Choosing $\phi:=v-u$, subtracting \eqref{q17} from the previous equation and applying H\"older's inequality, we get
\begin{align}
& \min\{\alpha_i,\alpha_j\}\int_{B_{\frac{r}{2}}}|\D(u-v)|^2\,dx\leq \max_{k=1\dots N}{\alpha_k}\Bigg[\int_{B_{\frac{r}{2}}\cap\mathcal{E}(j)\cap H}|\D u|^2\,dx+\int_{B_{\frac{r}{2}}\cap\mathcal{E}(i)\setminus H}|\D u|^2\,dx\\
& + \sum_{k\neq i,j}\bigg(\int_{B_{\frac{r}{2}}\cap \mathcal{E}(k)\cap H}|\D u|^2\,dx+\int_{B_{\frac{r}{2}}\cap \mathcal{E}(k)\setminus H}|\D u|^2\,dx+\int_{B_{\frac{r}{2}}\cap \mathcal{E}(k)}|\D u|^2\,dx\bigg) \Bigg]^{\frac{1}{2}}\bigg(\int_{B_{\frac{r}{2}}}|\D (u-v)|^2\,dx\bigg)^{\frac{1}{2}}.
\end{align}
Thus
\begin{align}
\int_{B_{\tau r}}|\D(u-v)|^2\,dx\leq C(\alpha)\bigg(\int_{B_{\frac{r}{2}}\cap(\mathcal{E}(j)\cap H)}|\D u|^2\,dx+\int_{B_{\frac{r}{2}}\cap(\mathcal{E}(i)\setminus H)}|\D u|^2\,dx+\sum_{k\neq i,j}\int_{B_{\frac{r}{2}}\cap\mathcal{E}(k)}|\D u|^2\,dx\bigg).
\end{align}
Arguing as above, the higher integrability of $\nabla u$ and H\"older's inequality yield
\begin{equation}
\int_{B_{\frac{r}{2}}\cap (\mathcal{E}(j)\cap H)}|\D u|^2\,dx\leq C^{\frac{1}{p}}\Bigg(\frac{|(\mathcal{E}(j)\cap H)\cap B_{\frac{r}{2}}|}{|B_{\frac{r}{2}}|}\Bigg)^{1-\frac{1}{p}}\int_{B_r}|\D u|^2\,dx,
\end{equation}
\begin{equation}
\int_{B_{\frac{r}{2}}\cap (\mathcal{E}(i)\setminus H)}|\D u|^2\,dx\leq C^{\frac{1}{p}}\Bigg(\frac{|(\mathcal{E}(i)\setminus H)\cap B_{\frac{r}{2}}|}{|B_{\frac{r}{2}}|}\Bigg)^{1-\frac{1}{p}}\int_{B_r}|\D u|^2\,dx,
\end{equation}
\begin{equation*}
\int_{B_{\frac{r}{2}}\cap \mathcal{E}(k)}|\D u|^2\,dx\leq C^{\frac{1}{p}}\bigg(\frac{|\mathcal{E}(k)\cap B_{\frac{r}{2}}|}{|B_{\frac{r}{2}}|}\bigg)^{1-\frac{1}{p}}\int_{B_r}|\D u|^2\,dx,\quad\forall k\neq i,j.
\end{equation*}
Choosing $\varepsilon_1$ such that $N\varepsilon_1^{1-\frac{1}{p}}=\tau^n$ and making the same computations as before we obtain the thesis.
\end{proof}

\begin{Thm}[Density upper bound]\label{Upper Bound}
Let $(\mathcal{E},u)$ be a $\Lambda$-minimizer of the functional ${\mathcal F}$ defined in \eqref{MainF}. Then, for any open set $U\Subset\Omega$ there exists a positive constant $C_2=C_2(n,N,\Lambda,\alpha)$ such that
\begin{equation}\label{decayF}
\mathcal{F}(\mathcal{E},u;B_r(x_0))\leq C_2 r^{n-1},
\end{equation}
for any $B_r(x_0)\subset U$.
\end{Thm}

\begin{proof}
Let $B_r(x_0)\subset U\Subset\Omega$. Let $i\in\{1,\dots,N\}$ such that $\alpha_i:=\min_{j\in\{1,\dots,N\}}\alpha_j$. We may assume that $x_0=0$ and $i=1$. We define:
\begin{align*}
& \mathcal{F}(1):=\mathcal{E}(1)\cup B_r,\\
& \mathcal{F}(h):=\mathcal{E}(h)\setminus B_r,\quad\forall h\in\{2,\dots,N\}.
\end{align*}
It holds that
\begin{align*}
& P(\mathcal{F}(1);\Omega)=P(\mathcal{E}(1)\setminus\overline{B_r};\Omega)+\mathcal{H}^{n-1}(\dd B_r\setminus\mathcal{E}(1)),\\
& \sum_{h=2}^N P(\mathcal{F}(h);\Omega)=\sum_{h=2}^N[P(\mathcal{E}(h)\setminus \overline{B_r};\Omega)+\mathcal{H}^{n-1}(\dd B_r\cap\mathcal{E}(h))]=\sum_{h=2}^N P(\mathcal{E}(h)\setminus \overline{B_r};\Omega)+\mathcal{H}^{n-1}(\dd B_r\setminus\mathcal{E}(1)).
\end{align*}
Using the previous equalities, the minimality of $(\mathcal{E},u)$ with respect to $(\mathcal{F},u)$ yields
\begin{align*}
\sum_{k=1}^N\alpha_k\int_{B_r\cap \mathcal{E}(k)}|\D u|^2\,dx+\frac{1}{2}\sum_{k=1}^N P(\mathcal{E}(k);B_r)\leq \alpha_1\int_{B_r}|\D u|^2\,dx+\mathcal{H}^{n-1}(\dd B_r\setminus\mathcal{E}(1))+2\Lambda |B_r|.
\end{align*}
It follows that
\begin{align*}\delta\alpha 
\int_{B_r\setminus\mathcal{E}(1)}|\D u|^2\,dx+\frac{1}{2}\sum_{k=1}^N P(\mathcal{E}(k);B_r)\leq c(n,\Lambda) r^{n-1}.
\end{align*}
where $\delta\alpha=\min_{k=2,\dots N}(\alpha_k-\alpha_1)$. We deduce that
\begin{equation}\label{Outof1}
\int_{B_r\setminus\mathcal{E}(1)}|\D u|^2\,dx+\sum_{k=1}^N P(\mathcal{E}(k);B_r)\leq c(n,\Lambda,\alpha) r^{n-1}.
\end{equation}
This concludes the proof of \eqref{decayF} in ${B_r\setminus\mathcal{E}(1)}$.\\
\indent It remains to prove that the Dirichlet integral decays in the right way also in ${B_r\cap\mathcal{E}(1)}$. This can be proved by contradiction using a quite standard blow-up argument that we present for the reader's convenience.\\
\indent We prove that there exist $\tau\in\big(0,\frac{1}{2}\big)$ and $M>0$ such that for every $\delta\in(0,1)$ there exists $h_0\in\N$ such that, for any $B_r(x_0)\subset U$, we have
\begin{equation}\label{contraddictionD}
\int_{B_r(x_0)}|\D u|^2\leq h_0r^{n-1} \quad\text{or}\quad \int_{B_{\tau r}(x_0)}|\D u|^2\,dx\leq M\tau^{n-\delta}\int_{B_r(x_0)}|\D u|^2\,dx.
\end{equation}
Arguing by contradiction, for every $\tau\in\big(0,\frac{1}{2}\big)$ and $M>0$ there exists $\delta\in(0,1)$ such that for every $h\in\N$ there exists a ball $B_{r_h}(x_h)\subset U$ such that
\begin{equation}
\label{a34}
\int_{B_{r_h}(x_h)}|\D u|^2\,dx>hr_h^{n-1}
\end{equation}
and
\begin{equation}
\label{a35}
\int_{B_{\tau r_h}(x_h)}|\D u|^2\,dx> M\tau^{n-\delta}\int_{B_{r_h}(x_h)}|\D u|^2\,dx.
\end{equation}
We choose $M>1$. Note that estimates \eqref{Outof1} and \eqref{a34} yield
\begin{equation}
\label{a44}
\sum_{k=2}^N\int_{B_{r_h}(x_h)\cap \mathcal{E}(k)}|\D u|^2\,dx+P(\mathcal{E}(k);B_{r_h}(x_h))\leq c_0r_h^{n-1}<\frac{c_0}{h}\int_{B_{r_h}(x_h)}|\D u|^2\,dx,
\end{equation}
and so
\begin{equation}
\label{a41}
\sum_{k=2}^N\int_{B_{r_h}(x_h)\cap \mathcal{E}(k)}|\D u|^2\,dx<\frac{c_0}{h}\int_{B_{r_h}(x_h)}|\D u|^2\,dx,
\end{equation}
for some positive constant $c_0$.\\
\indent Now we employ a typical blow-up argument. We set
\begin{equation*}
\varsigma_h^2:=\dashint_{B_{r_h}(x_h)}|\D u|^2\,dx
\end{equation*}
and, for $y\in B_1$, we introduce the sequence of rescaled functions defined as
\begin{equation*}
v_h(y):=\frac{u(x_h+r_hy)-a_h}{\varsigma_h r_h},\quad\text{with}\quad a_h:=\dashint_{B_{r_h}(x_h)}u\,dx.
\end{equation*}
We have ${\D u(x_h+r_hy)}=\varsigma_h\D v_h(y)$ and a change of variable yields
\begin{equation*}
\dashint_{B_1}|\D v_h(y)|^2\,dy=\frac{1}{\varsigma_h^2}\dashint_{B_{r_h}(x_h)}|\D u(x)|^2\,dx=1.
\end{equation*}
Therefore, there exist a (not relabeled) subsequence of $\left\{v_h\right\}_{h\in\N}$ and $v\in H^1(B_1)$ such that $v_h\rightharpoonup v$ in $H^1(B_1)$ and $v_h\rightarrow v$ in $L^2(B_1)$.
Moreover, the semicontinuity of the norm implies
\begin{equation}
\label{a39}
\dashint_{B_1}|\D v(y)|^2\,dy\leq\liminf_{h\rightarrow\infty}\dashint_{B_1}|\D v_h(y)|^2\,dy=1.
\end{equation}
Let us define the sets
\begin{equation*}
\mathcal{E}^*_h(k):=\frac{\mathcal{E}(k)-x_h}{r_h}\cap B_1.
\end{equation*}
We rewrite the inequalities \eqref{a34}, \eqref{a35} and \eqref{a41}. They become, respectively,
\begin{equation}
\label{a37}
\varsigma_h^2>\frac{h}{r_h},
\end{equation}
\begin{equation}
\label{a38}
\dashint_{B_\tau}|\D v_h(y)|^2\,dy>M\tau^{-\delta},
\end{equation}
\begin{equation}
\label{a42}
\sum_{k=2}^N\int_{B_1\cap \mathcal{E}^*(k)}|\D v_h(y)|^2\,dy<\frac{c_0}{h}\int_{B_1}|\D v_h(y)|^2\,dy=\frac{c_0\omega_n}{h}.
\end{equation}
To achieve the desired contradiction, it remains to show that the sequence $v_h$ cannot fulfill \eqref{a38} e \eqref{a42} because its limit $v$ minimizes the Dirichlet functional. Nevertheless, to establish a connection between the decay properties of $\nabla v$ and $\nabla v_h$
we must prove that the $L^2$-norm of $v_h$ converges to the $L^2$-norm of $v$. Observe that \eqref{a38} implies that $\varsigma_h\rightarrow\infty$, as $h\rightarrow\infty$. \\
\indent Since $r_h^{n-1} P(\mathcal{E}^*_h(k);B_1)=P(\mathcal{E}(k);B_{r_h}(x_h))$, by \eqref{a44}, we have that the sequence $\{P(\mathcal{E}^*_h(k);B_1)\}_{h\in\N}$ is bounded for any $k\in\{2,\dots,N\}$. Therefore, up to a not relabeled subsequence, $\mathbbm{1}_{\mathcal{E}^*_h(k)}\rightarrow\mathbbm{1}_{\mathcal{E}^*(k)}$ in $L^1(B_1)$, for some set $\mathcal{E}^*(k)\subset B_1$ of locally finite perimeter. By semicontinuity we deduce that
\begin{align}
\label{a43}
\int_{B_{1}}\mathbbm{1}_{\mathcal{E}^*(k)}|\nabla v|^2\,dy
& \leq \liminf_{h\rightarrow \infty}\int_{B_{1}}\mathbbm{1}_{\mathcal{E}^*(k)}|\nabla v_h|^2\,dy\\
& \leq \liminf_{h\rightarrow \infty}\Bigl(\int_{B_{1}}\mathbbm{1}_{\mathcal{E}_h^*(k)}|\nabla v_h|^2dy+\int_{B_{1}}\mathbbm{1}_{\mathcal{E}^*(k)\setminus \mathcal{E}_h^*(k) }|\nabla v_h|^2dy\Bigr)=0,
\end{align}
for any $k\in\{2,\dots,N\}$, where we used \eqref{a42} and the equi-integrability of $\big(|\nabla v_h|^2\big)_{h\in\N}$.
\\
\indent By $\Lambda$-minimality of $(\mathcal{E},u)$ with respect to $(\mathcal{E},u+\phi)$ we get, for $\phi\in H^1_0(B_{r_h}(x_h))$,
\begin{equation}
\sum_{k=1}^N\alpha_k\int_{B_{r_h}(x_h)\cap \mathcal{E}(k)}|\D u|^2\,dx\leq \sum_{k=1}^N\alpha_k\int_{B_{r_h}(x_h)\cap \mathcal{E}(k)}|\D u+\D \phi|^2\,dx   
\end{equation}
Using the change of variable $x=x_h+r_hy$, we deduce for every $\psi\in H^1_0(B_1)$
\begin{equation}
\sum_{k=1}^N\alpha_k\int_{B_1\cap \mathcal{E}_h^*(k)}|\varsigma_h\D v_h|^2\,dx\leq \sum_{k=1}^N\alpha_k\int_{B_1\cap \mathcal{E}_h^*(k)}|\varsigma_h\D v_h+\D \psi|^2\,dx.  
\end{equation}
Let $\eta\in C^\infty_c(B_1)$ such that $0\leq\eta\leq 1$. Choosing as a test function $\psi_h=\varsigma_h\eta(v-v_h)$,  we deduce that
\begin{align}
&\int_{B_1\cap \mathcal{E}_h^*(1)}\eta|\D v_h|^2\,dx\leq \int_{B_1\cap \mathcal{E}_h^*(1)}\eta|\D v|^2\,dx \\&+C(\alpha)\Biggl[\sum_{k=2}^N\int_{B_1\cap \mathcal{E}_h^*(k)}\Bigl(|\D v_h|^2+|\D v|^2\Bigr)\,dx +\int_{B_1}|\nabla \eta||v-v_h|\ dx\Biggr].
\end{align}
By \eqref{a42}, \eqref{a43} and the strong convergence of $v_h$ to $v$ in $L^2$ we deduce that the last term in the previous inequality tends to zero as $h\rightarrow\infty$, thus obtaining
\begin{equation}
\limsup_{h\rightarrow +\infty}{\int_{B_1}\eta|\D v_h|^2\,dx\leq \int_{B_1}\eta|\D v|^2\,dx }.    
\end{equation}
By the arbitrariness of $\eta$ and lower semicontinuity we conclude that $\D v_h\rightarrow\D v$ in $L^2$ and $v$ is harmonic. Thus, since $M>1$, by the harmonicity of $v$ and \eqref{a38} we get
\begin{equation*}
\dashint_{B_\tau}|\D v|^2\,dy\leq\dashint_{B_1}|\D v|^2\,dy\leq 1\leq \frac{\tau^\delta}{M}\lim_{h\rightarrow \infty}\dashint_{B_\tau}|\D v_h|^2\,dy<\lim_{h\rightarrow \infty}\dashint_{B_\tau}|\D v_h|^2\,dy=\dashint_{B_\tau}|\D v|^2\,dy,
\end{equation*}
which is a contradiction. Once \eqref{contraddictionD} is in force, the conclusion follows in a standard manner by applying an iteration argument.

\end{proof}
We are now in position to prove that in the neighborhoods where the perimeters of the chambers are small, the energy $\mathcal{F}$ decays in an appropriate manner.
\begin{Prop}
\label{DecadimEnergia}
Let $(\mathcal{E},u)$ be a $\Lambda$-minimizer of the functional ${\mathcal F}$ defined in \eqref{MainF}. Then, for any $\tau\in(0,1)$ there exists $\varepsilon_2=\varepsilon_2(\tau)>0$ such that, if for $B_r(x_0)\subset\Omega$
\begin{equation}
\label{q30}
P(\mathcal{E}(k);B_r(x_0))<\varepsilon_2 r^{n-1},\quad\forall k\in\{1,\dots,N\},
\end{equation}
then
\begin{equation}\label{decayE}
\mathcal{F}(\mathcal{E},u;B_{\tau r}(x_0))\leq C_3\tau^n(\mathcal{F}(\mathcal{E},u;B_r(x_0))+\Lambda r^n),
\end{equation}
for some positive constant $C_3=C_3(n,N,\Lambda,\alpha)$.
\end{Prop}

\begin{proof}
Let $B_r(x_0)\subset\Omega$ such that \eqref{q30} holds and $\tau\in(0,1)$. Setting
\begin{equation}
\mathcal{E}_r(k):=\frac{\mathcal{E}(k)-x_0}{r},\quad \mathcal{E}_r:=\{\mathcal{E}_r(k)\}_{k=1}^N,\quad u_r(y):=r^{-\frac{1}{2}}u(x_0+ry),\quad\forall y\in B_1,
\end{equation}
scaling \eqref{q30} we obtain
\begin{equation*}
P(\mathcal{E}_r(k);B_1)<\varepsilon_2,\quad\forall k\in\{1,\dots,N\}.
\end{equation*}
We need to prove that
\begin{equation*}
\mathcal{F}(\mathcal{E}_r,u_r;B_\tau)\leq\tau^n(\mathcal{F}(\mathcal{E}_r,u_r;B_1)+\Lambda r).
\end{equation*}
In what follows, for simplicity of notation, we will still denote $\mathcal{E}_r(k)$ by $\mathcal{E}(k)$ and $u_r$ by $u$ and we explicitly observe that $(\mathcal{E}(k),u)$ is a $\Lambda r$ minimizer of $\mathcal{F}$. \\
\indent We choose $\varepsilon_2<\min\Big\{\Big(\frac{\omega_n}{2c_{IP}N}\Big)^{\frac{n-1}{n}},\tau^{(n+1)(n-1)},\Big(\frac{\omega_n}{c_{IP}}\varepsilon_1\Big)^{\frac{n}{n-1}}\Big\}$, where $c_{IP}=c_{IP}(n)$ is the constant of the relative isoperimetric inequality and $\varepsilon_1=\varepsilon_1(\tau)$ is as in Lemma \ref{Lemma decadimento 1}. We first show that there exists $i\in\{1,\dots,N\}$ such that
\begin{equation}\label{q31}
|B_1\setminus\mathcal{E}(i)|\leq c_{IP}P(\mathcal{E}(i);B_1)^{\frac{n}{n-1}}.
\end{equation}
Let us assume by contradiction that the previous inequality is false for any $i\in\{1,\dots,N\}$. By the isoperimetric inequality, we have that
\begin{equation}
|\mathcal{E}(k)\cap B_1|\leq c_{IP}P(\mathcal{E}(k);B_1)^{\frac{n}{n-1}},\quad\forall k\in\{1,\dots,N\}.
\end{equation}
Thus, by the choice of $\varepsilon_2$, we get the following contradiction:
\begin{align*}
|B_1|=\sum_{k=1}^N|\mathcal{E}(k)\cap B_1|\leq c_{IP}\sum_{k=1}^N P(\mathcal{E}(k);B_1)^{\frac{n}{n-1}}\leq c_{IP}N\varepsilon_2^{\frac{n}{n-1}}<\frac{\omega_n}{2}.
\end{align*}
Therefore, \eqref{q31} is proved. We may assume that $i=1$ for simplicity. As a consequence we have:
\begin{equation}
\sum_{k=2}^N|\mathcal{E}(k)\cap B_1|=|B_1|-|\mathcal{E}(1)\cap B_1|=|B_1\setminus \mathcal{E}(1)|\leq c_{IP}P(\mathcal{E}(1);B_1)^{\frac{n}{n-1}}<\frac{\omega_n}{2}.
\end{equation}
The isoperimetric inequality yields
\begin{equation}
|\mathcal{E}(k)\cap B_1|\leq c_{IP}P(\mathcal{E}(k);B_1)^{\frac{n}{n-1}}, \quad\forall k\in\{2,\dots,N\}.
\end{equation}
Since
\begin{equation}
c_{IP}P(\mathcal{E}(1);B_1)^{\frac{n}{n-1}}\geq |B_1\setminus\mathcal{E}(1)|\geq \int_\tau^{2\tau}\mathcal{H}^{n-1}(\dd B_\rho\setminus \mathcal{E}(1))\,d\rho,
\end{equation}
we can choose $\rho\in(\tau,2\tau)$ such that
\begin{equation}
\label{q32}
\mathcal{H}^{n-1}(\dd B_\rho\setminus \mathcal{E}(1))\leq \frac{c_{IP}}{\tau}\varepsilon_2^{\frac{1}{n-1}}P(\mathcal{E}(1);B_1),
\end{equation}
\begin{equation}
\mathcal{H}^{n-1}(\dd^*\mathcal{E}(h)\cap \dd B_\rho)=0,\quad\forall h\in\{1,\dots,N\}.
\end{equation}
We set
\begin{align*}
& \mathcal{G}(1):=\mathcal{E}(1)\cup B_\rho,\\
& \mathcal{G}(k):=\mathcal{E}(h)\setminus B_\rho,\quad\forall k\in\{2,\dots,N\}.
\end{align*}
We remark that, $(\mathcal{G},u)$ is an admissible test pair to test the minimality of $(\mathcal{E},u)$ in $B_1$ because $\mathcal{E}(k)=\mathcal{G}(k)$ outside of $B_{\rho}$ for any $k\in\{1,\dots,N\}$. Thus
\begin{equation}\label{Mina}
\mathcal{F}(\mathcal{E},u;B_1)
 \leq \mathcal{F}(\mathcal{G},u;B_1)+ \Lambda r \rho^n.   
\end{equation}
To eliminate the common contribution in $B_1\setminus \overline{B}_{\rho}$ in the previous equation  we use the following equalities for the perimeter term:
\begin{align*}
P(\mathcal{E}(k);B_1)=P(\mathcal{E}(k);B_\rho)+P(\mathcal{E}(k);B_1\setminus\overline{B}_\rho),
\end{align*}
\begin{align*}
    P(\mathcal{G}(1);B_1)=P(\mathcal{E}(1);B_1\setminus \overline{B}_\rho)+\mathcal{H}^{n-1}(\dd B_\rho\setminus \mathcal{E}(1)),
\end{align*}
\begin{align*}
    \sum_{h=2}^N P(\mathcal{G}(h);B_1)
    & =\sum_{h=2}^N[P(\mathcal{E}(h);B_1\setminus \overline{B}_\rho)+\mathcal{H}^{n-1}(\dd B\rho\cap\mathcal{E}(h))]\\
    & =\sum_{h=2}^N P(\mathcal{E}(h);B_1\setminus \overline{B}_\rho)+\mathcal{H}^{n-1}(\dd B_\rho\setminus\mathcal{E}(1)).
\end{align*}
Moreover we observe that the choice of $\varepsilon_2$ implies that we are in position to apply Lemma \ref{Lemma decadimento 1}, since
\begin{equation*}
|\mathcal{E}(h)\cap B_1|\leq c_{IP}P(\mathcal{E}(h);B_1)^{\frac{n-1}{n}}\leq c_{IP}\varepsilon_2^{\frac{n-1}{n}}\leq\omega_n\varepsilon_1,\quad\forall h\in\{2,\dots,N\}.
\end{equation*}
Deleting the common term in \eqref{Mina} and applying Lemma \ref{Lemma decadimento 1} we conclude:
\begin{align*}
\mathcal{F}(\mathcal{E},u;B_\tau)
& \leq \mathcal{F}(\mathcal{E},u;B_\rho) \\  & \leq \alpha_1\int_{B_\rho}|\D u|^2\,dx+\mathcal{H}^{n-1}(\dd B_\rho\setminus\mathcal{E}(1))+c(n)\Lambda r\rho^n\\
& \leq \alpha_1\int_{B_{2\tau}}|\D u|^2\,dx+\frac{c_{IP}}{\tau}\varepsilon_2^{\frac{1}{n-1}}P(\mathcal{E}(1);B_1)+c(n)\Lambda r\tau^n\\
& \leq c(n,N,\alpha) \tau^n\int_{B_1}|\D u|^2\,dx+c(n)\tau^n P(\mathcal{E}(1);B_1)+c(n)\Lambda r\tau^n\\
& \leq c(n,N,\alpha)\tau^n(\mathcal{F}(\mathcal{E},u;B_1)+\Lambda r).
\end{align*}
\end{proof}
Building upon the previous proposition, we proceed to establish the following theorem, which states a lower bound estimate for the perimeter of the interfaces of optimal chambers.
\begin{Thm}[Density lower bound]
\label{Lower Bound}
Let $(\mathcal{E},u)$ be a $\Lambda$-minimizer of the functional ${\mathcal F}$ defined in \eqref{MainF} and $U\Subset\Omega$ be an open set. Then, there exists a positive constant $C_4=C_4(n,N,\Lambda,\alpha)$ such that, for every $x_0\in\bigcup_{k=1}^N\dd\mathcal{E}(k)\cap\Omega$ ad $B_r(x_0)\subset U$, it holds
\begin{equation}
\label{LowerBound}
\sum_{k=1}^N P(\mathcal{E}(k);B_r(x_0))\geq C_4 r^{n-1}.
\end{equation}
Moreover, $\mathcal{H}^{n-1}\Big( \Omega\cap\bigcup_{k=1}^N\dd\mathcal{E}(k)\setminus \bigcup_{k=1}^N\dd^*\mathcal{E}(k)\Big)=0$.
\end{Thm}
\begin{proof}
Since $\overline{\dd^*\mathcal{E}(k)}=\dd\mathcal{E}(k)$, for any $k\in\{1,\dots,N\}$, it is not restrictive to set $x_0=0\in\bigcup_{k=1}^N\dd^*\mathcal{E}(k)\cap\Omega$. 
Fix $\tau \in (0, 1)$ such that $2C_3 \tau^{1/2} < 1$, and fix $\sigma\in (0, 1)$ such that
\begin{equation}\label{small1}
2C_3 C_2 \sigma < \frac{\varepsilon_2(\tau)}{2},
\end{equation}
and let $r_0$ be such that
\begin{equation}\label{small2}
\Lambda r_0< \min\left\{ {\varepsilon_2(\tau)}, {C_2} \right\},
\end{equation}
where $C_2$, $C_3$ and $\varepsilon_2$ are the constants from Theorem \ref{Upper Bound} and Proposition \ref{DecadimEnergia}.
Assume by contradiction that for some $B_r \subset U$ with $r < r_0$, we have
\begin{equation}\label{small3}
\sum_{k=1}^N P(\mathcal{E}(k);B_r(x_0))\leq \varepsilon_2(\sigma) r^{n-1}.
\end{equation}
By induction, it is straightforward to show that
\begin{equation}\label{iter}
\mathcal{F}(\mathcal{E},u;B_{\sigma\tau^h r})\leq\varepsilon_2(\tau)\tau^{\frac{h}{2}}\big(\sigma\tau^h r\big)^{n-1},
\end{equation}
for every $h\geq 0$. Indeed, for $h=0$, using Proposition \ref{DecadimEnergia} and Theorem \ref{Upper Bound} and conditions \eqref{small1}, \eqref{small2} and \eqref{small3}
\begin{equation*}
  \mathcal{F}(\mathcal{E},u;B_{\sigma r})
  \leq C_3\sigma ^n(C_2 r^{n-1}+\Lambda r^{n})\leq 2C_3C_2 \sigma (\sigma r)^{n-1}\leq \varepsilon_2(\tau)(\sigma r)^{n-1}.
\end{equation*}
Assuming that \eqref{iter} holds for some $h$, to prove it also holds $h+1$ it suffices to apply Proposition \ref{DecadimEnergia} again and observe that $2C_2 \tau^{1/2} < 1$, together with \eqref{small2}. Indeed, we get
\begin{align}
\mathcal{F}(\mathcal{E},u;B_{\sigma\tau^{h+1} r})&\leq C_3 \tau^n \Big[\varepsilon_2(\tau)\tau^{\frac{h}{2}}\big(\sigma\tau^h r\big)^{n-1}+\Lambda(\sigma\tau^hr)^n\Big]\\
& \leq C_3 \varepsilon_2(\tau)\tau^n\Big[\tau^{\frac{h}{2}}\big(\sigma\tau^h r\big)^{n-1}+\sigma \tau^{h}\big(\sigma\tau^h r\big)^{n-1}\Big]\\
&\leq \varepsilon_2(\tau)2C_3\tau^{\frac 12}\tau^{\frac 12}\tau^{n-1}\tau^{\frac{h}{2}}\big(\sigma\tau^h r\big)^{n-1}\leq \varepsilon_2(\tau)\tau^{\frac{h+1}{2}}\big(\sigma\tau^{h+1} r\big)^{n-1}
\end{align}
It follows that
\begin{equation*}
\frac{1}{2}\sum_{k=1}^N P(\mathcal{E}(k);B_{\sigma\tau^h r})\leq \varepsilon_2(\tau)\tau^h\big(\sigma\tau^h r\big)^{n-1}.
\end{equation*}
Finally, it holds:
\begin{equation*}
\limsup_{\rho\rightarrow 0^+}\frac{\sum_{k=1}^N P(\mathcal{E}(k);B_{\rho})}{\rho^{n-1}}=\limsup_{h\rightarrow +\infty}\frac{\sum_{k=1}^N P(\mathcal{E}(k);B_{\sigma\tau^h r})}{\big(\sigma\tau^h r\big)^{n-1}}\leq 2\lim_{h\rightarrow+\infty}\varepsilon_2(\tau)\tau^h=0,
\end{equation*}
which is a contradiction.\\
\indent We are left to prove that 
\begin{equation}
\label{Equiv}
    \mathcal{H}^{n-1}\Big(\Omega\cap\bigcup_{k=1}^N\dd\mathcal{E}(k)\setminus \bigcup_{k=1}^N\dd^*\mathcal{E}(k)\Big)=0.
\end{equation}
By the lower bound \eqref{LowerBound}, we get
    \begin{align}
        &\limsup_{r\rightarrow 0^+}\frac{\mathcal{H}^{n-1}\big(\big(\bigcup_{k=1}^N\dd^*\mathcal{E}(k)\big) \cap B_r(x)\big)}{r^{n-1}}\\
        & =\limsup_{r\rightarrow 0^+}\sum_{k=1}^N\frac{\mathcal{H}^{n-1}(\dd^* \mathcal{E}(k)\cap B_r(x))}{r^{n-1}}=\limsup_{r\rightarrow 0^+}\sum_{k=1}^N\frac{P(\mathcal{E}(k);B_r(x))}{r^{n-1}}>0,\quad\forall x\in \bigcup_{k=1}^N\dd\mathcal{E}(k)\cap \Omega.
    \end{align}
On the other hand, by density property of $\mathcal{H}^{n-1}$-measurable sets with finite measure, see \cite[(2.42)]{AFP}, 
\begin{equation*}
 \limsup_{r\rightarrow 0^+}\frac{\mathcal{H}^{n-1}\big(\big(\bigcup_{k=1}^N\dd^*\mathcal{E}(k)\big) \cap B_r(x)\big)}{r^{n-1}} =0, \quad \text{for $\mathcal{H}^{n-1}$-a.e. }  x\notin \bigcup_{k=1}^N\dd^{*}\mathcal{E}(k)\cap \Omega,
\end{equation*}
thus \eqref{Equiv} follows.
\end{proof}

\subsection{Conclusion}

Theorem 1.1 summarizes the results obtained in Theorems 3.2 and 3.4, and is a direct consequence of them. We emphasize that the result concerns the union of the interfaces, but does not necessarily imply any regularity for the boundary of a single chamber. We believe that such a result could be achieved provided a suitable infiltration lemma were available. Specifically, if certain chambers of a minimizing partition occupy most of the ball 
$B_{2r}(x)$, then they must entirely fill the smaller ball $B_r(x)$ (cf. \cite[Lemma 30.2]{Ma} in the case of clusters). Establishing such a lemma, however, appears to be non-trivial due to the lack of sufficiently strong decay estimates for the Dirichlet bulk energy. However, its validity would represent a major breakthrough, potentially leading to more significant geometric insights into the chambers and, hopefully, to the regularity of the interface.\\

\emph{Acknowledgements.} The authors are members of the Gruppo Nazionale per l’Analisi Matematica, la Probabilità e le loro Applicazioni (GNAMPA) of the Istituto Nazionale di Alta Matematica (INdAM). The authors wish to acknowledge financial support from INdAM GNAMPA Project 2024 "Regolarità per problemi a frontiera libera e disuguaglianze funzionali in contesto finsleriano". The work of GP was partially supported by the project ”Harmony” (code MATE.IP.DR111.2024.HARMONY) within the program of the University ”Luigi Vanvitelli”.

\noindent
\author{Luca Esposito},
Dipartimento di Matematica, Università degli Studi di Salerno, Via Giovanni Paolo II 132, Fisciano 84084, Italy\\
luesposi@unisa.it

\medskip
\noindent
\author{Lorenzo Lamberti},
Institut \'Elie Cartan, Université de Lorraine, CNRS, IECL, F-54000 Nancy, France\\
lorenzo.lamberti@univ-lorraine.fr

\medskip
\noindent
\author{Giovanni Pisante}, Dipartimento di Matematica e Fisica, Universit\`{a} della Campania {\em Luigi Vanvitelli},
viale Lincoln 5, Caserta 81100, Italy\\
\noindent
giovanni.pisante@unicampania.it


\begin{thebibliography}{9}

\bibitem{ABr} L. Ambrosio and A. Braides. Functionals defined on partitions in sets of finite perimeter. II.
Semicontinuity, relaxation and homogenization. \emph{J. Math. Pures Appl.} {\bf 69} (1990), 307--333.

\bibitem{ABu} L. Ambrosio and G. Buttazzo, An optimal design problem with perimeter penalization, \emph{Calc. Var. Partial Differ. Equ.} \textbf{1} (1993), 55--69.

\bibitem{AFP} L. Ambrosio, N. Fusco and D. Pallara, \emph{Functions of Bounded Variation and Free Discontinuity Problems}, 1st ed., Oxford University Press, New York, 2000.

\bibitem{CEL1} {M. Carozza, L. Esposito, L. Lamberti}, Quasiconvex Bulk and Surface Energies:
$C^{1,\alpha}$ Regularity, \emph{Adv. Nonlinear Anal.} \textbf{13}(1) (2024), 20240021.

\bibitem{CEL2} M. Carozza, L. Esposito, L. Lamberti, Quasiconvex bulk and surface energies with subquadratic growth, \emph{Math. Eng.} \textbf{7}(3) (2025), 228--263.

\bibitem{CFP} M. Carozza, I. Fonseca and A. Passarelli Di Napoli, Regularity results for an optimal design problem with a volume constraint, \emph{ESAIM: COCV}, {\bf 20  no. 2} (2014), 460--487.

\bibitem{CFP2} M. Carozza, I. Fonseca and A. Passarelli Di Napoli, Regularity results for an optimal design problem with quasiconvex bulk energies, \emph{Calc. Var. Partial Differential Equ.}, {\bf 57}, 68 (2018).

\bibitem{CT} G. Congedo, I. Tamanini, On the existence of solutions to a problem in multidimensional segmentation. \emph {Ann. Inst. H. Poincaré Anal. Non Linéaire} {\bf 8} (1991), no. 2, pp. 175-–195 

\bibitem{DF} G. De Philippis and A. Figalli, A note on the dimension of the singular set in free interface problems, \emph{Differ. Integral Equ.} {\bf 28} (2015), 523--536.

\bibitem{DHV}{G. De Philippis}, {J. Hirsch} and G. Vescovo, {Regularity of minimizers for a model of charged droplets}, Commun. Math. Phys. \textbf{401}(1) (2023), 33--78.

\bibitem{E} L. Esposito, Density lower bound estimate for local minimizer of free interface problem with volume constraint, \emph{Ric. Mat.} \textbf{68, no. 2} (2019), 359--373.

\bibitem{EF} L. Esposito and N. Fusco, A remark on a free interface problem with volume constraint, \emph{J. Convex Anal.} {\bf 18 n.2} (2011), 417--426.

\bibitem{EL} L. Esposito and L. Lamberti, Regularity Results for an Optimal Design Problem with lower order terms, \emph{Adv. Calc. Var.}, \textbf{16}(4) (2023), 1093--1122.

\bibitem{EL2} {L. Esposito, L. Lamberti}, Regularity results for a free interface problem with H\"older coefficients, \emph{Calc. Var. Partial Diff. Equ.} {\bf 62} (2023), 156. https://doi.org/10.1007/s00526-023-02490-x

\bibitem{ELP} L. Esposito, L. Lamberti, G. Pisante, Epsilon-regularity for almost-minimizers of anisotropic free interface problem with H\"older dependence on the position. \emph{Interfaces Free Bound.} (2024), published online first. https://doi.org/10.4171/ifb/535
\bibitem{FJ} N. Fusco and V. Julin, On the regularity of critical and minimal sets of a free interface problem, \emph{Interfaces Free Bound.} \textbf{17 no.1} (2015), 117--142. 

\bibitem{H} T. Hales, The Honeycomb Conjecture. \emph{Discrete Comput. Geom.} {\bf 25} (2001), 1--22.

\bibitem{Lam} L. Lamberti, A regularity result for minimal configurations of a free interface problem, \emph{Boll. Un. Mat. Ital.} \textbf{14} (2021), 521-–539. 

\bibitem{LamLem} L. Lamberti and A. Lemenant, Quantitative regularity properties for the optimal design problem, preprint,
https://doi.org/10.48550/arXiv.2505.22365

\bibitem{Lar1} J. Larsen,  Distance between components in optimal design problems with
perimeter penalization, \emph{Ann. Sc. Norm. Super. Pisa, Cl. Sci., IV. Ser.} \textbf{28}(4) (1999), 641--649.

\bibitem{Lar2} C. J. Larsen, Regularity of components in optimal design problems with perimeter penalization, \emph{Calc. Var. Part. Diff. Equ.} {\bf 16} (2003), 17--29.

\bibitem{Leo} G. P. Leonardi, Infiltrations in immiscible fluids systems. \emph{Proceedings of the Royal Society of Edinburgh: Section A Mathematics.} {\bf 131(2)} (2001), 425--436.

\bibitem{Lin} F. H. Lin, Variational problems with free interfaces, \emph{Calc. Var. Part. Diff. Equ.} {\bf 1} (1993), 149--168.


\bibitem{LK} F. H. Lin and R. V. Kohn, Partial regularity for optimal design problems involving both bulk and surface
energies, \emph{Chin. Ann. of Math.} {\bf 20B} (1999), 137--158.

\bibitem{Ma} F. Maggi, \emph{Sets of finite perimeter and geometric variational problems. An introduction to geometric measure theory},  Cambridge Studies in Advanced Mathematics, 135. Cambridge University Press, Cambridge, 2012.

\bibitem{MV} E. Mukoseeva and G. Vescovo, Minimality of the ball for a model of charged liquid droplets, \emph{Ann. Inst. H. Poincaré Anal. Non Linéaire} {\bf 40} (2023), 457--509.

\bibitem{MN} Muratov, C.B., Novaga, M.: On well-posedness of variational models of charged drops. \emph{Proc. R. Soc. A} {\bf 472} (2016), 20150808.

\bibitem{TC} I. Tamanini and G. Congedo, Optimal segmentation of unbounded functions. \emph{Rend. Sem.
Mat. Univ. Padova} {\bf 95} (1996), 153--174.
2

\bibitem{W} B. White, Existence of least-energy configurations of immiscible fluids. \emph{J. Geom. Analysis}
{\bf 6} (1996), 151--161.
\end{thebibliography}
\end{document}